\documentclass[10pt]{amsart}
\usepackage{amscd}
\usepackage{amsfonts}
\usepackage{amssymb}
\usepackage{latexsym}

\newcommand{\ncm}{\newcommand}

\newtheorem{thm}{Theorem}[section]
\newtheorem{pro}[thm]{Proposition}
\newtheorem{lem}[thm]{Lemma}
\newtheorem{cor}[thm]{Corollary}
\newtheorem{lem&def}[thm]{Lemma \& Definition}
\newtheorem{defi}[thm]{Definition}
\newtheorem{exa}[thm]{Example}

\newtheorem{rmk}[thm]{Remark}

\def\M{\mathcal{M}}

\def\nn{\nonumber}

\ncm{\End}{\mbox{\rm End}\,}

\def\Hom{\mbox{\rm Hom}\,}

\def\id{\mbox{\rm id}}

\def\into{\hookrightarrow}
\def\to{\rightarrow}

\def\o{\otimes}    
\def\x{\times}     
\def\bra{\langle}
\def\ket{\rangle}

\ncm{\rarr}[1]{\stackrel{#1}{\longrightarrow}}
\ncm{\larr}[1]{\stackrel{#1}{\longleftarrow}}
\def\cop{\Delta}

\def\eps{\varepsilon}

\def\du1{\hat 1}
\def\Bra{\big\langle}
\def\Ket{\big\rangle}

\def\one{{\thinmuskip=5.5mu 1\!1\thinmuskip=3mu}}
\def\0{_{(0)}}
\def\1{_{(1)}}
\def\2{_{(2)}}
\def\3{_{(3)}}

\def\du1{\hat 1}
\def\lact{\triangleright}
\def\ract{\triangleleft}

\begin{document}

\title[Bialgebroids
and Depth Two Extensions]{Dual Bialgebroids for Depth Two Ring
Extensions}
\author{Lars Kadison}
\address{Matematiska Institutionen \\ G{\" o}teborg
University \\
S-412 96 G{\" o}teborg, Sweden}
\email{kadison@math.chalmers.se}
\author{Korn\'el Szlach\'anyi}
\address{Research Institute for Particle and Nuclear Physics, Budapest\\
H-1525 Budapest, P. O. Box 49, Hungary}
\email{szlach@rmki.kfki.hu}
\thanks{The first author thanks the second author and B. K\"ulshammer
for  pleasant visits to Budapest and Jena in the spring of 2001,
D. Nikshych, B. K\"ulshammer and O.A. Laudal for discussions,
and NorFA in Oslo for financial support.
The second author was partially supported by the Hungarian
Scientific Research Fund, OTKA T-034 512, and by NorFA in financing a visit
to Oslo.}
\subjclass{12F10,16W30, 22D30, 46L37}
\date{\today}

\begin{abstract}
We introduce a general notion of depth two for ring homomorphism $N \to M$,
and derive Morita equivalence of the step one and three centralizers,
$R = C_M(N)$ and $C = \End_{N-M}(M \o_N M)$, via dual bimodules
and step two centralizers
$A = \End\,_NM_N$ and $B = (M \o_N M)^N$, in a Jones tower above $N \to M$.
Lu's bialgebroids $\End_k A'$ and $A' \o_k {A'}^{\rm op}$ over a
$k$-algebra $A'$ are generalized to left and right bialgebroids
$A$ and $B$ with $B$ the $R$-dual bialgebroid of $A$.  We introduce
Galois-type actions of $A$ on $M$ and $B$ on $\End\,_NM$ when
$M_N$ is a balanced module.  In the case of Frobenius extensions $M | N$,
we prove an endomorphism ring theorem for depth two.  Further in
the case of irreducible extensions, we
extend previous results on Hopf algebra and weak
Hopf algebra actions in subfactor theory \cite{S,NV1}
and its generalizations \cite{KN,KN2} by methods other than nondegenerate
pairing. As a result, we have concrete
expressions for the Hopf or weak Hopf algebra structures on the step
two centralizers. Semisimplicity of $B$ is equivalent to
separability of the extension $M | N$.
In the presence of depth two, we show that biseparable
extensions are QF.
\end{abstract}
\maketitle

\section{Introduction}

Poisson and symplectic
groupoids were introduced by Weinstein in \cite{W1,W2} in the late eighties.
The notions extend to noncommutative algebra via Lu's notion of Hopf
algebroid
\cite{Lu} or bialgebroid with antipode. Some time before this,
Takeuchi \cite{T} introduced
the notion of $\times_R$-bialgebras based on studies of isomorphism classes
of simple algebras and earlier work by Sweedler \cite{Sw74}. A special case
of this extended notion of  bialgebra is Ravenel's commutative
Hopf algebroid introduced in the study of stable homotopy groups of spheres
\cite{R}.
Etingof and Varchenko \cite{EV} associated a Hopf algebroid to any dynamical
twist \cite{Ba}.

There is a quite different motivation coming from physics. In algebraic
quantum field theory \cite{Haag} the quest for finding a
2 dimensional analogue of the Doplicher-Roberts theorem \cite{DR}
(applies to quantum field theories in $d\geq3$ spacetime dimension) has lead
the authors of \cite{BSz} to introduce weak $C^*$-Hopf algebras (called also
quantum groupoids \cite{NV3}). The basic theory of weak Hopf algebras have
been
developed in  \cite{BNS,N,BSz2}. It turns out that bialgebroids and
$\times_R$-bialgebras are  equivalent\cite{X1,X2,BM}, while weak Hopf
algebras, Hayashi's face algebras, Maltsinotis's groupoid quantiques
and the finite dimensional version of Vallin's Hopf bimodules
occur as special cases \cite{EN, N, NV3, Vallin}.

In this paper we will bring the notion of bialgebroid together with a notion
of depth two in  the classification of subfactors \cite{P}.
Finite depth is a property of the standard
invariant of the Jones tower of  subfactor pair $N \subset M$ \cite{EK,JS}.
One forms the Jones
tower $N \subset M \subset M_1 \subset M_2 \subset \cdots$ by
iterating the basic construction $M_i = \bra M_{i-1}, e_i \ket$ where the
$e_i$ are the braidlike idempotents.  The tower of relative commutants
are the finite dimensional semisimple algebras $V_n = \Hom_{N-M_n}( M_n)$.
Finite depth is then the condition that the generating function $\sum_{n \geq
0}
\dim(V_n) x^n$ be rational; of depth  $n$
if $V_n = V_{n-1}e_n V_{n-1}$. Depth two inclusions are fundamental among
finite depth finite index inclusions via a Galois correspondence with
weak $C^*$-Hopf algebras and their coideal $*$-subalgebras \cite{NV2}
(in a role similar to Ocneanu's paragroups).

In this paper we investigate the fundamentals behind this trend beginning in
operator algebras and noncommutative Galois theory
of associating groups and quantum algebras to certain finite
depth algebra extensions.  Noncommutative
Galois theory for operator algebras was the
point of departure for Vaughan Jones's theory of subfactors
\cite{J}. In \cite{S,Longo} finite dimensional
Hopf $C^*$-algebras or Kac algebras are associated to
finite index irreducible subfactors of depth two. In \cite{NV1}
certain weak Hopf $C^*$-algebras \cite{BNS} are associated to non-integer
index subfactors
of finite depth with Galois correspondence \cite{NV2}.
In \cite{KN}
the depth two notion and results of \cite{S} are extended to an
algebraic analog without
trace: certain semisimple Hopf algebras are shown to have a Galois action
on split separable Frobenius extensions with trivial centralizer.
As we show in Section 8 of this paper,
a similar Hopf algebra $H$ may be associated
to an irreducible Frobenius extension $M | N$ of depth two;
semisimplicity or cosemisimplicity of $H$ being equivalent to
$M | N$ being a split or separable extension respectively. Even
the assumption of triviality in \cite{KN1}
for the centralizer may be relaxed
to separability (or absolute semisimplicity) at the price of obtaining
weak Hopf algebras, or quantum groupoids \cite{KN2}.
In each of these papers, it was essential to establish the quantum algebra
properties of $B$ together with its dual $A$ via a nondegenerate pairing.

In this paper we propose a completely general notion of depth two for
a ring extension $M | N$ which  allows the construction of bialgebroid
structures
on the centralizers $A$ and $B$ directly without a nondegenerate pairing.
In Section~2 we extend the theory of bialgebroids
to cover left and right bialgebroids and their duals, actions and smash
products.
We define depth~2 ring extension in Section~3 and derive from a certain
extension of Morita theory (cf.\ Section 1.1)
the basic classical properties among
the step~1, step~2, and step~3 centralizers in a Jones tower above a depth~2
extension $M | N$: the large centralizer $C$ is Morita equivalent to the
small centralizer $R$ (with no conditions imposed on it) while the
step~2 centralizers $A$ and $B$
are the Morita bimodules dual to one another and implementing the
equivalence.
In Section~4 we show
directly that $A$ is a left bialgebroid over $R$ with left action on
$M$:  if $M_N$ is balanced, the invariant subalgebra is $N$. We show that
$\End M_N$ is isomorphic to a smash product of $M$ with the bialgebroid $A$
over $R$, which is a basic step toward
a Galois theory for bialgebroid actions. In Section~5
we show directly that $B$ is a right bialgebroid  with
right action on $\End\,_NM$ and subalgebra of invariants $M$.
$A$ and $B$ are generalizations of Lu's  bialgebroids
in \cite[Section 3]{Lu} to noncommutative ring extensions,
and are shown to be $R$-dual to one another in Sections ~2,~3 and~5.
In Section~6 we specialize to the case where $M | N$ is
a Frobenius extension.  We answer a question in \cite{KN} by
showing that depth two passes up to the endomorphism
ring extension.    In Section~7 we show that
$A$ and $B$ specialize to isomorphic copies of the
dual Hopf algebras in \cite{KN} in case $R$ is
trivial in a depth two strongly separable extension of algebras.  We also
provide  an answer to a question
\begin{table}
\footnotesize
\begin{tabular}{|l|l|l|}      \hline \hline
\emph{Depth 2 Frobenius Extension} & \emph{Centralizer} & \emph{
Dual Quantum Algebras with Galois Actions} \\ \hline
biseparable algebra extension & trivial & semisimple, cosemisimple Hopf
algebra \\
algebra extension & trivial & Hopf algebra  \\
algebra extension  & separable & weak Hopf algebra \\
unrestricted & unrestricted & bialgebroid  \\ \hline \hline
\end{tabular}
\end{table}
in \cite{CK} in the presence of depth two by
showing that a  biseparable (i.e. split + separable + f.g. projective)
extension is quasi-Frobenius (QF). In Section~8 we extend
the results in \cite{KN} to an irreducible depth two
Frobenius extension by
finding an antipode $S : A \to A$ for a Frobenius bialgebroid over
trivial centralizer.  We prove that the action of
$A$ on $M$ in \cite{KN} is given by the analogous formula for the
action of $B$ on $M_1$.  In our last section, we generalize the main
results in \cite{KN2} by showing that $A$ and $B$ are weak Hopf algebras
dual to one another when $R$ is a separable algebra.
We summarize the algebraic results to date
in a table --- with a remark that there is in principle room for many more
entries in  future investigations.

\tableofcontents

In this paper rings are unital with $1 \neq 0$
 and  ring homomorphisms preserve the units.
A \textit{ring extension} $M|N$
in its most general sense is a ring homomorphism $\iota: N \to M$, which
induces a natural $N$-$N$-bimodule structure
on $M$ via $n \cdot m \cdot n' := \iota(n) m \iota(n')$;
the ring extension is
\textit{proper} if $N \into M$. The $\iota$
is suppressed in the language of ring extensions. $\,_M P_M$ denotes
an $M$-$M$-bimodule, and $P^M$ the centralizer subgroup $\{ p \in P|\,
pm = mp, \forall m \in M \}$, with the centralizer $C_M(N) = M^N$
being a special case.
A ring extension $M | N$ is said to have a property like (left) finitely
generated (f.g.) if ${}_NM$ and $M_N$ have this property (respectively, just
${}_NM$ is f.g.).  We denote $P$ being isomorphic to a direct
summand of another $M$-$M$-bimodule $Q$ by ${}_MP_M \oplus * \cong {}_MQ_M$.

\vskip 0.7truecm
\subsection{H-equivalence and Morita theory}

In this subsection, we recall a useful generalization of Morita theory which
will fit perfectly with the centralizer theory of depth two ring extensions
in Section~3.

Recall that two rings $T$ and $U$ are said to be Morita equivalent if
the category of left (or right) $T$-modules is equivalent to the category of
left (or right) $U$-modules.
The following statements due to Hirata \cite{Hi},
generalize Morita's main theorems in a very useful and simplifying
way.
Let $S$ be a ring with right modules $V$ and $W$.
\begin{lem}
\label{lem-hirata}
If $W_S \oplus * \cong (V \oplus \cdots \oplus V)_S$, and we set
$T := \End V_S$, $U := \End W_S$, then the natural bimodules ${}_TV_S$,
${}_UW_S$, ${}_U\Hom (V_S,W_S)_T$, ${}_T\Hom (W_S,V_S)_U$
are related by the following isomorphisms:
\begin{enumerate}
\item  a $U$-$U$-isomorphism
$\mu_T:  \Hom (V_S,W_S) \o_T \Hom (W_S,V_S) \stackrel{\cong}{\longrightarrow}
U$ via composition;
\item a $U$-$T$-isomorphism $\psi: \Hom (V_S,W_S)
\stackrel{\cong}{\longrightarrow} \Hom (\,_T\Hom(W_S,V_S),\,_TT)$
defined by $f \mapsto (g \mapsto g \circ f)$.
\item  a $U$-$S$-isomorphism $\tau: \Hom (V_S, W_S) \o_T V
\stackrel{\cong}{\longrightarrow} W$
given by $ f \o v \mapsto f(v)$.
\item a $U$-$S$-isomorphism $\iota: W \stackrel{\cong}{\longrightarrow} \Hom
(\,_T\Hom(W_S,V_S),\,_TV)$
via $w \mapsto (f \mapsto f(w))$.
\item $\Hom (V_S,W_S)$ is a finitely generated  projective right $T$-module
and a generator left $U$-module, while $\Hom (W_S, V_S)$ is a
finite projective  left $T$-module and generator right $U$-module.
\item $\Hom (\Hom (V_S,W_S)_T, \Hom (V_S,W_S)_T) \cong U$ and\hfill\break
      $\Hom (\,_U\Hom (W_S, V_S),\,_U\Hom (W_S, V_S))\cong U$.
\end{enumerate}
\end{lem}
\begin{proof} We observe that there are a finite number of
$f_i \in \Hom (V_S,W_S)$ and $g_i \in \Hom (W_S,V_S)$ such
that $\sum_i f_i \circ g_i = \id_W$. Define
$\mu_T^{-1}(u) = \sum_i u f_i \o g_i$.  Define
$\psi^{-1}(F) = \sum_i f_i \circ F(g_i)$.
The rest of the proof
is quite similar and left to the reader \cite{Hi}.
\end{proof}

The lemma has the following easy converse:  if $\mu_T$ is epi, then
$W_S \oplus * \cong \oplus^n V_S$.
The lemma leads directly
to Hirata's result \cite{Hi}:
\vskip 0.7truecm
\begin{pro}
If both $V \oplus * \cong \oplus^n W$ and $W \oplus * \cong \oplus^m V$
(in which case we say $V$ and $W$ are H-equivalent $S$-modules),
then $T$ and $U$ are Morita equivalent rings with Morita context
given by $$({}_U\Hom (V_S,W_S)_T, \,_T\Hom (W_S,V_S)_U, \mu_T, \mu_U).$$
\end{pro}

If $V_S$ is f.g. projective and a generator,
then $V_S$ and $S_S$ are H-equivalent, and $\End V_S$ is Morita equivalent
to $S$ via $V$ and its right $S$-dual $V^*$. This recovers Morita's theorem.

\section{Bialgebroids}

In this  section, we define left and right ring-theoretic versions of Lu's
$R$-bial\-ge\-broid$A$,
an $A$-module algebroid $M$ which has a measuring action and subring of
invariants, and
a smash product ring $M \rtimes A$. In the final subsection, we introduce
the $R$-dual right bialgebroids $B$ and $B'$ of a left bialgebroid $A$, whose
multiplicative and comultiplicative structures are defined on
the  right and left dual spaces $\Hom (A_R,R_R)$ and $\Hom ({}_RA,{}_RR)$.
This section covers technical preliminaries about bialgebroids that might
be consulted when needed.

Lu's original definition  of a bialgebroid \cite{Lu} corresponds to our
 left bialgebroid below if all maps are in the category of
$k$-algebras. The necessity to introduce a left and a
right version comes from the asymmetry of the bialgebroid axioms
under the switch
to the opposite ring structure. The axioms we use here for the right
bialgebroid are those of \cite{Sz}, easily seen to be equivalent
to the right-handed version of Lu's original axioms.

Let $R$ be a ring. A {\bf right bialgebroid} over $R$ consists of the data
and
axioms:
\renewcommand{\theenumi}{\roman{enumi}}
\begin{enumerate}
\item a ring $A$ and two ring homomorphisms $R^{op}\rarr{t}A\larr{s}R$
such that $s(r')t(r)$ $=t(r)s(r')$ for $r,r'\in R$. Thus $A$ can be made into
an $R$-$R$-bimodule by setting $r\cdot a\cdot r':=at(r)s(r')$.
\item $R$-$R$-bimodule maps $\cop\colon A\to A\o_R A$ and $\eps\colon A\to R$
such that the triple $\bra A,\cop,\eps\ket$ is a comonoid in the category
$_R\M_R$. (Another name is $R$-coring.)
\item $\cop$ is multiplicative in the following sense. Although $A\o_R A$ has
no ring structure in general, its sub-bimodule
$$
A\times_R A:=\{X\in A\o_R A\,|\,(s(r)\o 1)X=(1\o t(r))X,\,\forall r\in R\,\}
$$
is a ring with multiplication $(a\o a')(a''\o a''')=aa''\o a'a'''$. Now we
require that
$$
\cop\colon A\to A\times_R A
$$
be a ring homomorphism.
\item $\eps$ preserves the unit: $\eps(1)=1_R$
\item $\eps$ is compatible with the ring structure on $A$ in the sense of the
axioms
$$
\eps(t(\eps(a))b)\ =\ \eps(ab)\ =\ \eps(s(\eps(a))b)\,,\qquad a,b\in A\,.
$$
\end{enumerate}
When discussing duals of bialgebroids in Subsection~\ref{sec: duals} we shall
see that property (v) is dual to the unitalness of the coproduct, $\cop(1)=1\o
1$, which is part of property (iii) above.

Without much comment we  list the axioms of a {\bf left
bialgebroid} over $R$. It consists of
\begin{enumerate}
\item a pair $R\rarr{s}A\larr{t}R^{op}$ of ring homomorphisms such that
$s(r)t(r')=t(r')s(r)$, $r,r'\in R$
\item $R$-$R$-bimodule maps $\cop\colon A\to A\o_R A$ and $\eps\colon A\to R$
where $A$ is given the bimodule structure $r\cdot a\cdot r':=s(r)t(r')a$
\end{enumerate}
such that
\begin{enumerate}
\item[(ii/1)] $(\cop\o\id_A)\circ\cop=(\id_A\o\cop)\circ\cop$
\item[(ii/2)] $(\eps\o\id_A)\circ\cop=\id_A=(\id_A\o\eps)\circ\cop$
\item[(iii/1)] $\cop(a)(t(r)\o 1)=\cop(a)(1\o s(r))$, $a\in A$, $r\in R$
\item[(iii/2)] $\cop(ab)=\cop(a)\cop(b)$, $a,b\in A$
\item[(iii/3)] $\cop(1)=1\o 1$
\item[(iv)] $\eps(1)=1_R$
\item[(v)] $\eps(a s(\eps(b)))=\eps(ab)=\eps(a t(\eps(b)))$, $a,b\in A$.
\end{enumerate}
\renewcommand{\theenumi}{\arabic{enumi}}

It follows from the $R$-linear property of $\cop$ that
\begin{equation}\label{eq:del s}
\cop(s(r)) = s(r) \o 1,
\end{equation}
\begin{equation}\label{eq:del t}
\cop(t(r)) = 1 \o t(r),
\end{equation}
for a left
bialgebroid.
Note that our convention is that $s$ is a homomorphism and $t$ is an
anti-homomorphism from $R$ both for left and right bialgebroids.
In the language of weak Hopf algebras \cite{BNS} $s(R)$ corresponds to $A^L$
in the case of left bialgebroids but corresponds to $A^R$ in the case of right
bialgebroids. In particular, we see
the formulas (\ref{eq:del s}) and (\ref{eq:del t}) are
interchanged for a right bialgebroid.

As for the relation of left and right bialgebroids we note that if
$\bra A,R,s,t,\cop,\eps\ket$ is a left bialgebroid then
$\bra A',R,s',t',\cop,\eps\ket$ is a right bialgebroid where
\begin{equation}
A'\ =\ A^{op},\quad s'=t^{op}\colon R\to A^{op},\quad t'=s^{op}\colon
R^{op}\to
A^{op}.
\end{equation}

On the other hand, passing to the opposite coring structure does not change
"handedness". As a matter of fact if
$\bra A,R,s,t,\cop,\eps\ket$ is a left bialgebroid then $A^{cop}:=\bra
A,R',s',t',\cop',\eps\ket$ is also a left bialgebroid where
\begin{equation}
R'=R^{op},\quad s'=t,\quad t'=s
\end{equation}
thus the bimodule structure of $_{R'}{A'}_{R'}$ is the opposite of $_RA_R$,
i.e.,
\begin{equation}
r_1\cdot' a\cdot' r_2\ =\ t(r_1)s(r_2)a\ =\ r_2\cdot a\cdot r_1\,,\quad
r_1,r_2\in R,\ a\in A\,.
\end{equation}
Applying the Sweedler notation $a\1\o a\2$ for $\cop(a)$ the coproduct of
$A^{cop}$ is
\begin{equation}
\cop'=\cop^{op}\colon a\mapsto a\2\o_{R^{op}}a\1
\end{equation}
for which $\eps$ is the counit.

The same construction yields a right bialgebroid $A^{cop}$ from a right
bialgebroid $A$.

\begin{exa}
\begin{rm}
If $A$ is an algebra over a commutative ring $R$ with $s = t = u: R \to A$,
the unit map, then a left
or right $R$-bialgebroid structure on $A$ is just a bialgebra \cite{Mo}.
\end{rm}
\end{exa}

\begin{exa}
\begin{rm}
A weak Hopf algebra $A$ with left (or target) subalgebra $A^L$
(a separable algebra)
is a bialgebroid (indeed Hopf algebroid \cite{Lu}) over $A^L$ \cite{EN}.
Conversely, if $A$ is a bialgebroid over a separable $K$-algebra $R$, where
$K$ is a field, then $A$ has a weak Hopf $K$-algebra structure given
in Proposition~\ref{pro: WHA} and \cite{Sz}.
\end{rm}
\end{exa}

\subsection{Module algebroids}
The extra structure on a ring $A$ which makes it a left (right)
bialgebroid over $R$ is precisely a monoidal structure on its category
$_A\M$ ($\M_A$) of modules together with a strictly monoidal forgetful functor
to $_R\M_R$ (cf. \cite{Sz}). Therefore the natural candidate for a "module
algebra" over a bialgebroid $A$ is a monoid in the category of $A$-modules.
More explicitely a {\bf left $A$-module algebroid} over a left bialgebroid
$\bra A,R,s,t,\cop,\eps\ket$ consists of
\begin{itemize}
\item a left $A$-module $_AM$ inheriting an $R$-$R$ bimodule structure from
the $A$-action: $r\cdot m\cdot r'=(r\cdot 1_A\cdot r')\lact m=s(r)t(r')\lact
m$
\item an associative multiplication $\mu_M\colon M\o_R M\to M$, $m\o
m'\mapsto
mm'$ satisfying
\begin{equation} \label{mul M}
a\lact(mm')\ =\ (a\1\lact m)(a\2\lact m'),\quad a\in A,\ m,m'\in M
\end{equation}
\item and a unit $\eta_M\colon R\to M$, $r\mapsto r\cdot 1_M\equiv 1_M\cdot
r$
for the multiplication $\mu_M$ which satisfies
\begin{equation} \label{unit M}
a\lact 1_M\ =\ \eps(a)\cdot 1_M\,,\qquad a\in A\,.
\end{equation}
\end{itemize}

Note then that
\begin{equation}\label{*}
(m \cdot r)m' = m(r \cdot m'),
\end{equation}
as well as $r \cdot (mm') = (r \cdot m)m'$ and
$(mm') \cdot r = m(m' \cdot r)$.
Notice that Eqns.\ (\ref{mul M}) and (\ref{unit M}) express the fact that $\bra
M,\mu_M,\eta_M\ket$ is not only a monoid in $_R\M_R$ but also in $_A\M$.

A {\bf right $A$-module algebroid} over a right bialgebroid $\bra
A,R,s,t,\cop,\eps\ket$ consists of
\begin{itemize}
\item a right $A$-module $M_A$ inheriting an $R$-$R$ bimodule structure from
the $A$-action: $r\cdot m\cdot r'=m\ract(r\cdot 1_A\cdot r')=m\ract
t(r)s(r')$
\item an associative multiplication $\mu_M\colon M\o_R M\to M$, $m\o
m'\mapsto
mm'$ satisfying
\begin{equation}
(mm')\ract a\ =\ (m\ract a\1)(m'\ract a\2),\quad a\in A,\ m,m'\in M
\end{equation}
\item and a unit $\eta_M\colon R\to M$, $r\mapsto r\cdot 1_M\equiv 1_M\cdot
r$
for the multiplication $\mu_M$ which satisfies
\begin{equation}
1_M\ract a\ =\ 1_M\cdot\eps(a)\,,\qquad a\in A\,.
\end{equation}
\end{itemize}

A left $A$-module algebroid over the left bialgebroid $A$ is the same as the
right $A^{op}$-module algebroid over the right bialgebroid $A^{op}$.

If $_AM$ is a left $A$-module algebroid over the left bialgebroid $A$ then
the
opposite ring $M^{op}$ yields a monoid in $_{R^{op}}M_{R^{op}}$ such that
$M^{op}$ becomes a left $A^{cop}$-module algebroid.

Similarly, a right $A$-module algebroid $M_A$ gives rise to a right
$A^{cop}$-module algebroid $M^{op}_{A^{cop}}$.

\subsection{The subring of invariants}
Let $_AM$ be a module algebroid over the left bialgebroid $\bra
A,R,s,t,\cop,\eps\ket$. The invariants of $M$ is the subset
\begin{equation} \label{def: inv}
M^A:=\{n\in M|a\lact m=s(\eps(a))\lact m, a\in A\}.
\end{equation}
Notice that if $n\in M^A$ then
\begin{equation}
t(\eps(a))\lact n=s(\eps(t(\eps(a))))\lact n=s(\eps(a))\lact n=a\lact n
\end{equation}
for all $a\in A$. Therefore we obtain an equivalent definition if $s$ is
replaced by $t$ in (\ref{def: inv}).
\begin{lem}\label{lem: inv<End NMN}
For $m\in M$, $n\in M^A$ and $a\in A$ we have
\begin{equation}
a\lact(mn)=(a\lact m)n\,,\qquad a\lact(nm)=n(a\lact m)\,.
\end{equation}
In particular, it follows that $M^A$ is a subring of $M$.
\end{lem}
\begin{proof}
In the next calculation, we use one of the two equivalent definitions of
invariants, next the identity $m_1m_2=(1\1\lact m_1)(1\2\lact m_2)$,
then axiom
(iii) of a left bialgebroid and finally one of the counit axioms.
\begin{eqnarray*}
a\lact(mn)&=&(a\1\lact m)(s(\eps(a\2))\lact n)=(1\1 a\1\lact m)(1\2
s(\eps(a\2))\lact n)\\
&=&(1\1 t(\eps(a\2)) a\1\lact m)(1\2\lact n)=(a\1\cdot\eps(a\2)\lact m)n\\
&=&(a\lact m)n,
\end{eqnarray*}
and
\begin{eqnarray*}
a\lact(nm)&=&
t(\eps(a\1))\lact n)(a\2 \lact m) = (n \cdot \eps(a\1))(a\2 \lact m) \\
&=&n(s(\eps(a\1))a\2 \lact m)= n(a\lact m). \qed
\end{eqnarray*}
\renewcommand{\qed}{}\end{proof}

In a similar way the invariants of a right bialgebroid $M_A$ can be written
in
two ways
\begin{equation}
M^A=\{n\in M|m\ract a=m\ract s(\eps(a)),\ a\in A\}=\{n\in M|m\ract a=m\ract
t(\eps(a)),\ a\in A\}
\end{equation}
and form a subring $M^A\subset M$.

Another important subring in a module algebroid is the sub-$A$-module
generated by the identity. For a left $A$-module algebroid $M$ it is
\begin{equation}
M^j:=\{a\lact 1_M|a\in A\}.
\end{equation}
It is the image of the map
\begin{equation}
j_M\colon R\to M,\qquad r\mapsto s(r)\lact 1_M
\end{equation}
which is a ring homomorphism. As a matter of fact, $t(r)\lact
1_M=s(\eps(t(r)))\lact 1_M=s(r)\lact 1_M$ therefore
\begin{eqnarray*}
j_M(r)j_M(r')&=&(s(r)\lact 1_M)(t(r')\lact 1_M)\\
&=&((r\cdot 1\cdot r')\1\lact 1_M)((r\cdot 1\cdot r')\2\lact 1_M)\\
&=&s(r)t(r')\lact 1_M=s(r)s(r')\lact 1_M=s(rr')\lact 1_M\\
&=&j_M(rr')\,.
\end{eqnarray*}
As a consequence of the lemma above, $M^j$ commutes with the invariants,
\begin{equation}
(s(r)\lact 1_M)n=n(s(r)\lact 1_M)\,\quad r\in R,\ n\in M^A\,.
\end{equation}
For the module algebroids we shall consider in Sections \ref{sec: A} and
\ref{sec: B} the $M^j$ is  actually equal to the centralizer of $M^A$ in $M$.

\subsection{The smash product}
If $A$ is a left bialgebroid over $R$ and $M$ is a left $A$-module algebroid
then $M$ is a right $R$ module via $m\cdot r:=m j_M(r)$.

\begin{defi}
The smash product $M\rtimes A$ of a left $A$-module algebroid ${}_AM$ with
$A$
is the ring with additive group  $M\o_R A$ and multiplication
defined by
\begin{equation}
(m\rtimes a)(m'\rtimes a'):=m(a\1\lact m')\rtimes a\2a'\,.
\end{equation}
Analogously one defines $A\ltimes M$ for a right $A$-module algebroid $M_A$.
\end{defi}

The multiplication is well-defined because of Eq.\ (\ref{*}) and (2.iii).
The maps $\imath_M\colon m\mapsto m\rtimes 1$ and $\imath_A\colon a\mapsto
1_M\rtimes a$ are ring homomorphisms of $M$, respectively of $A$, into
$M\rtimes A$. One can check easily the relations
\begin{eqnarray}
\imath_M(m)\imath_A(a)&=&m\rtimes a\,, \label{eq: imath} \\
\imath_A(a)\imath_M(m)&=&(a\1\lact m)\rtimes a\2\,.
\end{eqnarray}
for $m\in M$, $a\in A$. The $\imath_M$ is always an embedding by the following
argument. Lemma \ref{lem: inv<End NMN} allows us to map $M\rtimes A$ into
$\End M_N$ where $N:=M^A$, the subring of invariants, with $A$ mapped into
$\End(\,_NM_N)\subset \End M_N$. As a matter of fact $m\rtimes a$ acts on $M$
as $\lambda(m)(a\lact -)$. Composing the ring map $M\rtimes A\to \End
M_N$ with $\imath_M$ one obtains $\lambda_M$, the left regular representation
of $M$, which is faithful. Thus $\imath_M$ must be mono and the smash product
is always a proper ring extension of $M$.

On the other hand, $\imath_A$ is not necessarily mono.
If $\imath_A(a)=0$ then using the above map into $\End\,_NM_N$ again
we obtain that $a\lact m=0$ for all $m\in M$. By
Eq.\ (\ref{eq: imath}), $m\rtimes a=0$ for all
$m\in M$. So  if either $_AM$ is faithful or if $M_R$ is
faithfully flat, then $A$ embeds into the smash product $M\rtimes A$ via
$\imath_A$.

\subsection{Duals}\label{sec: duals}

If $A$ is a bialgebroid over $R$ one may expect a bialgebroid structure on
the
dual bimodule $A^*$ or $^*A$ provided $A_R$ or $_RA$ is finitely projective.
The fine point here is that in taking duals one really has to take into
account that $A$ is not only a bimodule over $R$ but carries 4 actions of
$R$:
multiplying either from the left or right by either $s_A(r)$ or $t_A(r)$.
Comultiplications of left and right bialgebroids are bimodule maps with
respect to two different (and disjoint) pairs of $R$-actions. Multiplication,
however, cannot be written as a bimodule map in either of these two
categories but requires the use of "mixed" pairs of $R$-actions. This is why
in defining duals of a bialgebroid we have to use new bimodule structures
of $A$ and not those appearing before in comultiplications.

\subsubsection{The right dual $A^*$}
Let $A$ be a left bialgebroid over $R$ and assume that $A_R$ is finitely
generated projective. Recall that the right action of $r\in R$ is $a\mapsto
t_A(r)a$. We shall denote by $A^{(t)}$ the $R$-$R$-bimodule which is the
additive group $A$ on which $r\in R$ acts from the left via $a\mapsto
at_A(r)$
and acts from the right via $a\mapsto t_A(r)a$. Thus the
right $R$-action of $A^{(t)}$ coincides with the right $R$-action on $A$
dictated by the left bialgebroid structure. But the left action is
different. We define the right dual of $A$ as the
right dual bimodule of $A^{(t)}$, i.e., $A^*=\Hom(A_R,R_R)$ carrying the
bimodule structure
\begin{equation}
\bra r\cdot b\cdot r',a\ket:=r\bra b,at_A(r')\ket\,,\qquad b\in A^*,\ a\in
A\,.
\end{equation}
Here and below $\bra b,a\ket$ denotes the canonical pairing, i.e., the
evaluation of $b$ on $a$. Now we make $A^*$ into a ring by defining
multiplication via the formula
\begin{equation}\label{eq: three}
\bra bb',a\ket\ :=\ \bra b',\bra b,a\1\ket\cdot a\2\ket\,.
\end{equation}
which is associative due to coassociativity of $\cop_A$.
Note with caution that $\cdot$ here denotes the ordinary $R$-bimodule
structure on $A$: $r \cdot a = s(r)a$.
The multiplication has a unit $1_{A^*}=\eps_A$.

If $A^*$ is going to be a right bialgebroid then the maps
\begin{eqnarray}
s_{A^*}\ \colon R\to A^*\,,&\quad&r\mapsto 1_{A^*}\cdot
r=\eps_A\circ\rho_A(t_A(r)) \label{eq: two}\\
t_{A^*}\,\colon R^{op}\to A^*\,,&\quad& r\mapsto r\cdot
1_{A^*}=\eps_A\circ\lambda_A(s_A(r))
\end{eqnarray}
are ring homomorphisms. That this is indeed the case follows from
previous identities such as Eq. (\ref{eq:del t}):
\begin{eqnarray}
\bra bs_{A^*}(r),a\ket&=& \bra s_{A^*}(r),\bra b,a\1\ket\cdot a\2\ket
                       =\eps_A(\bra b,a\1\ket\cdot a\2t_A(r))\nn\\
&=&\bra b,a\1\cdot\eps_A(a\2t_A(r))\ket=\bra b,at_A(r)\ket\ =\ \bra b\cdot
                                                             r,a\ket\\
\bra bt_{A^*}(r),a\ket&=&\bra t_{A^*}(r),\bra b,a\1\ket\cdot a\2\ket
                       =r\eps_A(\bra b,a\1\ket\cdot a\2)\nn\\
&=&r\bra b,a\1\ket\eps_A(a\2)=r\bra b,a\ket\ =\ \bra r\cdot b,a\ket
\end{eqnarray}
For future convenience we list the five basic symmetry relations of the
pairing, two of which have just been proved:
\begin{eqnarray}
\bra b,t_A(r)a\ket&=&\bra b,a\ket r\\
\bra b,s_A(r)a\ket&=&\bra b,\eps_A(s_A(r)a\1)\cdot a\2\ket=\bra b,\bra
t_{A^*}(r),a\1\ket\cdot a\2\ket\nn\\
&=&\bra t_{A^*}(r)b,a\ket\\
\bra b,at_A(r)\ket&=&\bra bs_{A^*}(r),a\ket\\
\bra b,as_A(r)\ket&=&\bra b,\eps_A(a\1)\cdot a\2s_A(r)\ket=\bra b,\eps_A(a\1
t_A(r))\cdot a\2\ket\nn\\
&=&\bra s_{A^*}(r)b,a\ket\\
\bra bt_{A^*}(r),a\ket&=&r\bra b,a\ket
\end{eqnarray}
In order to define comultiplication on $A^*$ we have to utilize that $A_R$ is
finitely projective. A consequence of this is that the natural map
\begin{eqnarray}
A^*\o_R A^*&\to&\Hom((A^{(t)}\o_R A^{(t)})_R,R_R)\\
b\o b'&\mapsto&\{a'\o a\mapsto\bra b\cdot\bra b',a'\ket,a\ket\}\nn
\end{eqnarray}
is an isomorphism. Its inverse can be given in terms of dual bases $\{a_i\}$
of $A_R$ and $\{b_i\}$ of $_RA^*$ as $f\mapsto\sum_{i,j}f(a_j\o a_i)\cdot
b_i\o b_j$.

Noticing that for any $b\in A^*$ the map $a'\o a\mapsto\bra b,aa'\ket$
belongs
to the above hom-group, a comultiplication
\begin{equation}
\cop_{A^*}\colon A^*\to A^*\o_R A^*\,,\qquad b\mapsto b\1\o b\2
\end{equation}
can be defined by requiring
\begin{equation}\label{eq: one}
\bra b\1\cdot\bra b\2,a'\ket,a\ket\ =\ \bra b,aa'\ket\,,\qquad b\in A^*,\
a,a'\in A\,.
\end{equation}
In terms of the dual bases it can be written as
\begin{equation}
\cop_{A^*}(b)\ =\ \sum_{i,j}\ \bra b,a_ia_j\ket\cdot b_i\o b_j\ .
\end{equation}
Now we turn to verifying the bialgebroid axioms for $\cop_{A^*}$.

$\cop_{A^*}$ is a bimodule map by its very definition
(a simple calculation with Eq.\ \ref{eq: one}). In order to see that
its image lies in $A^*\x_R A^*$ we compute
\begin{eqnarray}
\bra s_{A^*}(r)b\1\cdot\bra b\2,a'\ket,a\ket&=&\bra b\1\cdot\bra
b\2,a'\ket,as_A(r)\ket\nn\\
&=&\bra b,as_A(r)a'\ket=\bra b\1\cdot \bra b\2,s_A(r)a'\ket,a\ket\nn\\
&=&\bra b\1\cdot\bra t_{A^*}(r)b\2,a'\ket,a\ket
\end{eqnarray}
Now it is meaningful to ask whether the map $\cop_{A^*}\colon A^*\to A^*\x_R
A^*$ is a ring homomorphism. The proof of multiplicativity goes as follows
(lines 6 to 7 below requires Eqs.\ \ref{eq: two}, \ref{eq: three}
and~\ref{eq:del t}):
\begin{eqnarray*}
\bra (bb')\1\cdot\bra (bb')\2,a'\ket,a\ket&=&\bra bb',aa'\ket
=\bra b',\bra b,(aa')\1\ket\cdot(aa')\2\ket\\
&=&\bra b',\bra b,a\1a'\1\ket\cdot a\2a'\2\ket\\
&=&\Bra b',\left(\bra b\1\cdot
\bra b\2,a'\1\ket,a\1\ket\cdot a\2\right)a'\2\Ket\\
&=&\Bra {b'}\1\cdot\bra {b'}\2,a'\2\ket,\bra b\1\cdot\bra
b\2,a'\1\ket,a\1\ket\cdot
a\2\Ket\\
&=&\Bra\left(b\1\cdot\bra b\2,a'\1\ket\right)\left({b'}\1\cdot\bra
{b'}\2,a'\2\ket\right),a\Ket\\
&=&\Bra b\1 s_{A^*}(\bra b\2,a'\1\ket){b'}\1 s_{A^*}(\bra
{b'}\2,a'\2\ket),a\Ket\\
&=&\Bra b\1 {b'}\1 s_{A^*}\left(\Bra t_{A^*}(\bra b\2,a'\1\ket){b'}\2,a'\2\Ket
\right),a\Ket\\
&=&\bra b\1 {b'}\1\cdot\bra {b'}\2,s_A(\bra b\2,a'\1\ket)a'\2\ket,a\ket\\
&=&\bra b\1 {b'}\1\cdot\bra {b'}\2,\bra b\2,a'\1\ket\cdot a'\2\ket,a\ket\\
&=&\bra b\1 {b'}\1\cdot \bra b\2 {b'}\2,a'\ket,a\ket
\end{eqnarray*}
Preservation of the unit, $\cop_{A^*}(1)=1\o 1$, can be seen as
\begin{eqnarray*}
\bra 1\1\cdot\bra 1\2,a'\ket,a\ket&=&\bra 1,aa'\ket=\eps_A(aa')=
\eps_A(at_A(\eps_A(a')))\\
&=&\bra 1,at_A(\bra 1,a'\ket)\ket=\bra 1s_{A^*}(\bra 1,a'\ket),a\ket\\
&=&\bra 1\cdot\bra 1,a'\ket,a\ket
\end{eqnarray*}
We are left with constructing the counit for $A^*$. Let
\begin{equation}
\eps_{A^*}\colon A^*\to R\,,\qquad b\mapsto\bra b,1_A\ket\,.
\end{equation}
Then $\eps_{A^*}$ is an $R$-$R$-bimodule map,
\begin{eqnarray*}
\eps_{A^*}(r\cdot b\cdot r')&=&\bra bs_{A^*}(r')t_{A^*}(r),1_A\ket
=r\bra b,t_A(r')\ket\\
&=&r\eps_{A^*}(b)r'\,,
\end{eqnarray*}
it satisfies the counit properties because
\begin{eqnarray*}
\bra b\1\cdot\eps_{A^*}(b\2),a\ket&=&\bra b\1\cdot\bra b\2,1_A\ket,a\ket
=\bra b,a\ket\\
\bra \eps_{A^*}(b\1)\cdot b\2,a\ket&=&\bra b\1,1_A\ket\bra b\2,a\ket
=\bra b\1,t_A(\bra b\2,a\ket)\ket\\
&=&\bra b\1\cdot\bra b\2,a\ket,1_A\ket=\bra b,a\ket\,,
\end{eqnarray*}
it preserves the unit,
\[
\eps_{A^*}(1_{A^*})\ =\ \eps_A(1_A)\ =\ 1_R\,,
\]
and finally it is compatible with multiplication of $A^*$,
\begin{eqnarray*}
\eps_{A^*}(s_{A^*}(\eps_{A^*}(b))b')&=&\bra s_{A^*}(\bra
b,1_A\ket)b',1_A\ket\\
&=&\bra b',s_A(\bra b,1_A\ket)\ket=\bra b',\bra b,1_A\ket\cdot 1_A\ket
=\bra bb',1_A\ket\\
&=&\eps_{A^*}(bb')\\
\eps_{A^*}(t_{A^*}(\eps_{A^*}(b))b')&=&\bra b',s_A(\bra
b,1_A\ket)\ket \\
&=&\eps_{A^*}(bb')\,.
\end{eqnarray*}

What we have just proven is the following:

\begin{pro}\label{pro: rightdual}
If $A$ is a left bialgebroid over $R$ such that $A_R$ is finitely generated
projective then $B:=\Hom(A_R,R_R)$ has a unique right bialgebroid structure
over $R$ such that
\begin{eqnarray}
\bra bb',a\ket&=&\bra b',\bra b,a\1\ket\cdot a\2\ket\\
\bra b,aa'\ket&=&\bra b\1\cdot\bra b\2,a'\ket,a\ket
\end{eqnarray}
where $\bra\ ,\ \ket\colon B\x A\to R$ denotes the canonical pairing.
\end{pro}

\subsubsection{The left dual $^*A$}
Let $A$ again be a left bialgebroid but now assume that $_RA$ is finitely
generated projective. For $a\in A$ and $b\in \,^*A=\Hom(\,_RA,\,_RR)$ we
denote by $[a,b]\in R$ the evaluation of $b$ on $a$. As a bimodule $^*A$ is
considered to be the dual bimodule of $A^{(s)}$ where the latter is the
additive group $A$ on which $r\in R$ acts from the left by $a\mapsto s_A(r)a$
and from the right by $a\mapsto as_A(r)$. Then similarly as in the above
Proposition we can construct a right bialgebroid structure on $^*A$. More
precisely we have

\begin{pro}\label{pro: leftdual}
If $A$ is a left bialgebroid over $R$ such that $_RA$ is finitely generated
projective then $^*A:=\Hom(\,_RA,\,_RR)$ has a unique right bialgebroid
structure over $R$ such that
\begin{eqnarray}
[a,bb']&=&[a\1\cdot[a\2,b],b']\\
{[}aa',b]&=&[a,[a',b\1]\cdot b\2]\,.
\end{eqnarray}
\end{pro}

The proof is very similar to the previous construction and therefore omitted.
We only give here the symmetry properties of the $[\ ,\ ]$ pairing:
\begin{eqnarray}
[s_A(r)a,b]&=&r[a, b] \\
{[}t_A(r)a,b]&=&[a,bt_{^*A}(r)]\\
{[}as_A(r),b]&=&[a,s_{^*A}(r)b]\\
{[}at_A(r),b]&=&[a,bs_{^*A}(r)]\\
{[}a,t_{^*A}(r)b]&=&[a,b]r
\end{eqnarray}

We note that both the $\bra\ ,\ket$ and $[\ ,\ ]$ pairings are variations of
Schauenburg's skew pairing $\tau$ of \cite{Sch} with the caution that
\cite{Sch} uses only left bialgebroids in our language.

\subsubsection{Duals of right bialgebroids}
Left and right duals $^*B$ and $B^*$ of a right bialgebroid $B$ can be
introduced directly using the above notions of duals of left bialgebroids. Let
$B$ be a right bialgebroid over $R$ with $B_R$ finitely generated projective.
Then its right dual $B^*$ is a left bialgebroid defined by
\begin{equation}
B^*:=\left((B^{\rm op})^*\right)^{\rm op}
\end{equation}
This means that $B^*=\Hom(B_R,R_R)$ as an additive group and its bialgebroid
structure is to be read from the canonical pairing
$$
[a,b]:=a(b)\ \mbox{for}\ a\in B^*,\ b\in B
$$
satisfying precisely the relations of the pairing $[\ ,\ ]$ of Proposition
\ref{pro: leftdual}. Now it is easy to verify that for a left bialgebroid $A$
such that $_RA$ is finitely generated projective the canonical isomorphism
$A\cong (^*A)^*$ of Abelian groups is in fact an isomorphism of left
bialgebroids. In other words, if $B$ is the left dual of $A$ then $A$ is the
right dual of $B$.

The same conclusion holds for a left bialgebroid $A$ such that $A_R$ is finite
projective. Its left dual $B=\,^*A$ is a right bialgebroid such that $_RB$ is
finite projective and for such a $B$ a left dual can be introduced via
\begin{equation}
^*B:=\left(\,^*(B^{\rm op})\right)^{\rm op}\,.
\end{equation}
Denoting the canonical pairing of $b\in B$ with $a\in\, ^*B$ by $\bra b,a\ket$
we obtain that $\bra\ ,\ \ket$ satisfies the relations of Proposition
\ref{pro: rightdual}. Thus again if $B$ is the left dual of the left
bialgebroid $A$ then $A$ is the right dual of $B$.

In Sections \ref{sec: A} and \ref{sec: B} we shall meet a situation when
the left bialgebroid $A$ has both a left and a right dual and they are
isomorphic to a right bialgebroid $B$. In this case it is fair to say that $A$
and $B$ are dual pairs of bialgebroids.

\section{Depth 2 Ring Extensions}

In this pivotal section, we extend the notion of depth two from subfactors
\cite{GHJ}and Frobenius extensions \cite{KN} to a conciser notion for any
ringextension or homomorphism $N \to M$. Depth two has equivalent
formulations in terms of H-equivalence and quasibasis.
We note that the  H-separable
extension defined in \cite{Hi} is a particular example --- and has a
certain parallel  theory developed by Sugano. The ``step two centralizers''
$A := \End {}_NM_N$ and $B := (M \o_N M)^N$ play a large role in the
theory of depth two extensions.  Our main theorem shows that the centralizer
$R := M^N$ and ``step three'' centralizer $\End {}_N(M \o_N M)_M$ are
Morita equivalent with the step two centralizers  as the Morita
bimodules, and provides several results  needed in later sections.

The tensor-square $M \o_N M$, left and right endomorphism rings,
$\End\,_NM$ and $\End M_N$, of a ring extension $M|N$
have the natural $M$-$M$ bimodule structures given by
$m \cdot x \o y \cdot m' = mx \o ym'$, $(m\eta m')(x) = \eta(xm)m'$
\newline
and $(m f m')(x) = m f(m'x)$ for $m,m',x,y \in M$, $\eta \in \End\,_NM$,
and $f \in  \End M_N$, respectively.

\begin{defi}\label{defi: D2}
A ring extension $M|N$ is called left depth two or
{\bf left D2} if $$_NM\o_N M_M\oplus * \cong \oplus^n\,_NM_M$$
for some positive integer $n$;
{\bf right D2} if $$_MM\o_N M_N\oplus * \cong \oplus^m\,_MM_N$$ for some
positive integer $m$.

$M|N$ is called {\bf D2} if it is D2 both from the left and from the
right.
\end{defi}

In particular, the natural modules ${}_M M \o M$ and $M \o M_M$ are f.g.
projective
for a D2 extension $M | N$.

\begin{rmk}
\begin{rm}
If left D2,
$M$ and $M \otimes_N M$ are in fact H-equivalent as $M \o N^{\rm
op}$-modules,
since the multiplication mapping $\mu_N: M \o_N M \to M$ is a split
$N$-$M$-epi
for any ring extension $M | N$.  A similar statement is equivalent to the
right D2 condition.
\end{rm}
\end{rmk}
\begin{exa}
\begin{rm}
A classical depth two subfactor $N \subseteq M$
 of finite index is of depth two
by Proposition~\ref{pro-classical D2}.
\end{rm}
\end{exa}
\begin{exa}
\label{exa-cpre}
\begin{rm}
A centrally projective ring extension $M|N$ is D2,
since $\,_N M _N \oplus * \cong \oplus^n \,_NN_N$
and we may arrive at the definition above by tensoring from the left or
right by $\,_MM_N$ or $\,_NM_M$.
In particular, if $N \subset Z(M)$, the center of $M$, and
$M$ is  finitely
generated and projective as a $N$-module,  then $M$ is a D2 extension of $N$,
a ``D2 algebra'' over $N$.
\end{rm}
\end{exa}

\begin{exa}
\begin{rm}
An H-separable extension  \cite{Hi} $M|N$ is D2, since its defining property
is
$$_M(M\o_N M)_M\oplus * \cong \oplus^n\,_MM_M.$$
\end{rm}
\end{exa}

Let $A := \End \,_NM_N$ and $B := (M \o_N M)^N$.

\begin{lem}\label{lem: D2}
$N\subset M$ is left D2 iff there exist $b_i\in B$ and $\beta_i\in A$
(called a left D2 quasibasis) such that
$$
\sum_ib^1_i\o b_i^2\beta_i(m)\ =\ m\o 1\,,\qquad m\in M\,.
$$
$N\subset M$ is right D2 iff there exist $c_i\in B$ and $\gamma_i\in A$ such
that
$$
\sum_i\gamma_i(m)c_i^1\o c_i^2\ =\ 1\o m\,,\qquad m\in M\,.
$$
\end{lem}
\begin{proof}
(Left D2 $\Rightarrow$ existence of quasibasis.)
Let $\pi: \oplus^n\,_NM_M \to \,_NM \o M_M$ and $\sigma: \,_NM \o M_M \to
\oplus^n\,_NM_M$
denote the split epi and its section implied by the definition. Furthermore,
let
$\{ e_i \}_{i=1}^n$ be the standard basis of the free module $\oplus^n M_M$,
and $p_i : \oplus^n\,_NM_M \to \,_NM_M$ be the standard projections.
Then we let $b_i := \pi(e_i)$ where clearly $b_i \in B$.  If
$\iota_1: \,_NM_N \into \,_NM \o_N M_N$ denotes the map $m \mapsto m \o 1$,
then we let $\beta_i := p_i \circ \sigma \circ \iota_1 \in A$.
Then
$$ m \o 1 = \pi (\sigma (\iota_1(m))) = \sum_i \pi(e_i) \beta_i(m) = \sum_i
b_i^1 \o b_i^2 \beta_i(m)
$$
The rest of the proof is similar.
\end{proof}

\begin{rmk}
\begin{rm}
If $M$ is a finite dimensional $k$-algebra with $N =k1$, with dual bases $e_i
\in M$
and $\pi_i \in M^* = \Hom_k(M,k)$, then a left and right D2 quasibasis is
given by $\beta_i = \gamma_i = \iota \pi_i$, $b_i= e_i \o 1_M$ and $c_i = 1_M
\o e_i$,
where $\iota: k \into M$ is the unit map.
\end{rm}
\end{rmk}

For every $m \in M$, we let $\lambda(m) \in \End M_N$ denote $\lambda(m)(x) =
mx$
and $\rho(m) \in \End\,_NM$ denote $\rho(m)(x) = xm$.  If $r \in R :=
C_M(N)$,
we note that $\lambda(r), \rho(r) \in \End\,_NM_N = (\End\,_NM)^N = (\End
M_N)^N$.
In the sequel the $R$-bimodule structure on $A$ is understood to be $r \cdot
\alpha \cdot r' =  \lambda(r) \rho(r') \alpha$.

\begin{pro}\label{pro: end}
If $M|N$ is a right D2 extension, then $\End\,_NM \cong A \o_R M$ as
$N$-$M$-bimodules via
$\alpha \o m \mapsto \rho(m) \alpha$.  If $M | N$ is a left D2 extension,
then
$\End M_N \cong M \o_R A$ as $M$-$N$-bimodules via $m \o \alpha \mapsto
\lambda(m) \circ \alpha$.
\end{pro}
\begin{proof}
We claim that $f \mapsto \sum_i \gamma_i \o c_i^1 f(c_i^2)$
for $f \in \mathcal{E}' := \End\,_NM$ defines an inverse.
Since $\sum_i \gamma_i(m) c_i^1 f(c_i^2) = f(m)$, we see that
$f = \sum_i \rho(c_i^1 f(c_i^2)) \gamma_i \in \rho(M)A$.

Similarly an inverse to the second statement is given by $$f \mapsto
\sum_i f(b^1_i) b^2_i \o \beta_i$$
for each $f \in \mathcal{E} := \End M_N$.
\end{proof}

\begin{pro}\label{pro: tensor-square}
If $M|N$ is left or right D2, then
\begin{equation}
A \otimes_R A \cong \Hom_{N-N}(M \o_N M,M)
\end{equation}
via $\alpha \o \beta \longmapsto ( m \o m' \mapsto \alpha(m) \beta(m'))$.
\end{pro}
\begin{proof}
  The inverse mapping is given by
\begin{equation}
\Hom_{N-N}(M \o_N M,M) \to A \otimes_R A,\ \
f \longmapsto \sum_i f(- \o b_i^1)b^2_i \o_R \beta_i,
\end{equation}
since
$$\sum_i \alpha(-) \beta(b_i^1)b_i^2 \o \beta_i = \sum_i \alpha \o
\beta(b_i^1)b_i^2 \beta_i(-) = \alpha \o \beta
$$
and $\sum_i f(m \o b_i^1) b_i^2  \beta_i(m') = f(m \o m')$ for each $m,m' \in
M$.
We can carry out a similar proof with a right D2 quasibasis.
\end{proof}

Next is a main theorem for depth two extensions.   We make use
of the ``step one'' centralizer $R$, and ``step two centralizers'' $A$
and $B$ defined above; in addition, a ``step three'' centralizer
$C := \End\,_N(M \o_N M)_M$ (cf.\ \cite{GHJ}, \cite{KN}).

\begin{thm}\label{thm: D2}
If $M | N$ is left D2, then $C$ and $R$ are Morita equivalent rings
with invertible bimodules ${}_CB_R$ and ${}_RA_C$ in a Morita context.
In particular, $B_R$ and $\,_RA$ are f.g. projective generators
with the following isomorphisms:
\begin{enumerate}
\item $\mu_R: B \o_R A \stackrel{\cong}{\longrightarrow} C$ via
$b \o \alpha \mapsto ( m \o m' \mapsto b \alpha(m)m')$.
\item $\psi: B_R \stackrel{\cong}{\longrightarrow} \Hom (\,_RA,\,_RR)_R$
via $b \mapsto (\alpha \mapsto \alpha(b^1)b^2)$.
\item $\tau: B \o_R M \stackrel{\cong}{\longrightarrow} M \o_N M$
defined by $\tau(b \o m) = bm$.
\item $\iota: M \o_N M \stackrel{\cong}{\longrightarrow} \Hom (\,_RA,\,_RM)$
via $\iota(m \o m')(\alpha) = \alpha(m)m'$.
\item $C \cong \End B_R$ via $c \mapsto (b \mapsto c(b))$.
\item $C \cong \End\,_RA$ via
\begin{equation}\label{eq: nn}
c \longmapsto (\alpha \mapsto \mu(\alpha \o \id_M)c\iota_1)
\end{equation}
where $\iota_1 : M \to M \o_N M$ by $m \mapsto m \o 1$.
\end{enumerate}
\end{thm}
\begin{proof}
First we note that $R \cong \End_{N-M}(M)$ via $r \mapsto \lambda(r)$
with inverse $f \mapsto f(1)$.  Next we
note that $\Hom_{N-M} (M, M \o_N M) \cong B$ via $f \mapsto f(1)$
with inverse $b \longmapsto (m \mapsto bm).$  The bimodule structure
on ${}_CB_R$ is given by
$$ c \cdot b \cdot r = c(b^1 \o b^2 r).  $$

We next note that $A \cong \Hom_{N-M} (M \o_N M, M)$ via $\alpha \mapsto (m \o
m' \mapsto \alpha(m)m')$ with inverse $f \mapsto f \circ \iota_1$.   The
bimodule structure on ${}_RA_C$ is given by
$$ r \cdot \alpha \cdot c = \lambda(r) \mu_N(\alpha \o \id_M) c\iota_1\,.$$
The rest follows strictly from the Lemma and Proposition in
the introduction;
however, we note some useful inverses to some of the isomorphisms above.
\begin{equation}
\tau^{-1}(m \o m') =  \sum_i b_i \o \beta_i(m)m'
\end{equation}
\begin{equation}
\iota^{-1}(f) =  \sum_i b^1_i \o b^2_i f(\beta_i)
\end{equation}
\begin{equation}
\psi^{-1}(\phi) =  \sum_i b_i \phi(\beta_i)
\end{equation}
Dual bases for $\,_RA$ are given by $\{ \psi(b_i) \}$, $\{ \beta_i \}$.
\end{proof}

By yet another application of Lemma~\ref{lem-hirata}
we  prove in a similar way
(but writing arguments to the left of a function) that
if $M | N$ is \textit{right} D2, then the natural module  ${}_RB$
and $A_R$, where $\alpha \cdot r = \rho(r) \circ \alpha$,
are progenerators with corresponding isomorphisms,
such as $M \o_R B \cong M \o M$ via $m \o b \mapsto mb$,
\begin{equation}
\label{B dual to A}
B \cong \Hom(A_R,R_R), \ \ b \longmapsto (\alpha \mapsto b^1 \alpha(b^2))
\end{equation}
and
\begin{equation}
C \cong \End A_R.
\end{equation}

 From Prop.~\ref{pro: end} and the theorem we easily establish
\begin{cor}
\label{cor-endo}
If $M | N$ is D2, then
\begin{equation}
{}_N\mathcal{E}'_M \oplus * \cong \oplus^s {}_NM_M,
\end{equation}
\begin{equation}
{}_M\mathcal{E}_N \oplus * \cong \oplus^t {}_MM_N.
\end{equation}
\end{cor}

We obtain a type of converse to the theorem by
noting that if $\mu_R$ is epi, then $M | N$ is left D2. Equivalently, ${}_RA$
f.g. projective,
$\psi$ an isomorphism with $C \cong \End\,_RA$ implies $M | N$ is left D2.
This shows that a classical depth two pair of semisimple algebras is
left D2 and similarly right D2 \cite{GHJ}.

\section{The Left Bialgebroid $A$ and its Action}\label{sec: A}

In this section, we construct a left bialgebroid structure on the ring
$A=\End {}_NM_N$ given a D2 ring extension $M|N$. This bialgebroid
recovers Lu's endomorphism algebra example if $N$ is trivial.
Its underlying coproduct is given in terms of a left or right D2 quasibasis
defined in the previous section. Its action on $M$ is the natural action
of endomorphisms. The fixed points of this action will be the image of $N$
in $M$if $M_N$ is balanced, whose definition is recalled below. The right
endomorphism ring is in either case isomorphic to the smash product ring
$M\rtimes A$. For the next theorem, we recall that a right
$R$-module $V$ is \textit{balanced}
if the natural left $\mathcal{E} := \End V_R$-module on $V$
has left endomorphism ring naturally
anti-isomorphic to $R$:  $R \stackrel{\cong}{\to} \End_{\mathcal{E}} V$.
For example, if $V_R$ is a generator (i.e.
$R_R \oplus * \cong \oplus^n V_R$), then $V_R$ is balanced by
the well-known Morita's lemma.

\begin{thm} \label{thm - A}
Let $N \to M$ be a depth two extension of rings.
Then  $A$ is a left bialgebroid over $R$ with left action of $A$ on
$M$.  If $M_N$ is moreover balanced,
then the subring of invariants under this action is $N$.

More explicitly, the bialgebroid is $\bra A, R, s_A, t_A, \cop_A, \eps_A\ket$
where
\begin{eqnarray}
A&=&\End\,_NM_N\\
R&=&C_M(N)\\
s_A(r)&=&\lambda(r)\colon m\mapsto rm\\ \label{eq: alf}
t_A(r)&=&\rho(r)\colon m\mapsto mr\\
r\cdot\alpha\cdot r'&=&\lambda(r)\rho(r')\alpha\colon
m\mapsto r\alpha(m)r'\\
\cop_A(\alpha)&=&\sum_i \gamma_i\o_R c_i^1\alpha(c_i^2 -)\\
\eps_A(\alpha)&=&\alpha(1_M)
\end{eqnarray}
The $A$-module action on $M$ is simply the action of endomorphisms, $\alpha
\lact m=\alpha(m)$.
\end{thm}

\begin{proof}At first we check the left bialgebroid axioms.

$\cop_A$ is an $R$-$R$-bimodule map:
\begin{eqnarray*}
\cop_A(r\cdot \alpha\cdot r')&=&\gamma_j\o_R c_j^1r\alpha(c_j^2 -)r'\\
&=&\gamma_j\o_R c_j^1r\alpha(c_j^2b_i^1)b_i^2\beta_i(-)r'\\
&=&\gamma_j(-)c_j^1r\alpha(c_j^2b_i^1)b_i^2\o_R \beta_i(-)r'\\
&=&r\cdot\left(\alpha(- b_i^1)b_i^2\o_R \beta_i(-)\right)\cdot r'
\end{eqnarray*}
Putting $r=r'=1$ yields an alternative formula for the coproduct,
\begin{equation}\label{eq: copA}
\cop_A(\alpha)=\alpha(- b_i^1)b_i^2\o_R \beta_i
\end{equation}
which, when plugged back, gives
\[
\cop_A(r\cdot \alpha\cdot r')=r\cdot\cop_A(\alpha)\cdot r'
\]

Coassociativity:
\begin{eqnarray*}
(\cop_A\o\id_A)\circ\cop_A(\alpha)&=&\gamma_j\o_R c_j^1\alpha(c_j^2-
b_i^1)b_i^2\o_R\beta_i\\
(\id_A\o\cop_A)\circ\cop_A(\alpha)&=&\gamma_j\o_R c_j^1\alpha(c_j^2 -
b_i^1)b_i^2\o_R\beta_i
\end{eqnarray*}

The property $\cop_A(\alpha)(\rho(r)\o 1)=\cop_A(\alpha)(1\o\lambda(r))$:
\begin{eqnarray*}
LHS&=&\alpha(- rb_i^1)b_i^2\o_R \beta_i=\\
&=&\alpha(- b_j^1)b_j^2\beta_j(rb_i^1)b_i^2\o_R \beta_i\\
&=&\alpha(- b_i^1)b_i^2\o_R\beta_i(r-)\\
&=&RHS
\end{eqnarray*}

Multiplicativity of $\cop_A$:
\begin{eqnarray*}
\cop_A(\alpha)\cop_A(\alpha')&=&
\alpha(\alpha'(- b_j^1)b_j^2b_i^1)b_i^2\o_R\beta_i(\beta_j(-))\\
&=&\gamma_k\o_R
c_k^1\alpha(\alpha'(c_k^2b_j^1)b_j^2b_i^1)b_i^2\cdot\beta_i(\beta_j(-))\\
&=&\gamma_k\o_R c_k^1\alpha(\alpha'(c_k^2-))\\
&=&\cop_A(\alpha\circ\alpha')
\end{eqnarray*}

Unitalness: $\cop_A(1)=1\o_R 1$ and $\eps_A(1)=1_R$ are obvious.

The compatibility of $\eps_A$ with multiplication:

\begin{equation*}
\eps_A(\alpha\circ\lambda(\eps_A(\alpha')))=\alpha(\alpha'(1_M))=
\eps_A(\alpha\alpha')
\end{equation*}
and the same for $\rho$ instead of $\lambda$.

This completes the proof that $A$ is a bialgebroid.

Module algebra properties:
$A$ acts on $M$ by the simple formula $\alpha\triangleright m:=\alpha(m)$.
The induced $R$-$R$-bimodule structure on $M$ is also the obvious one arising
 from $R$ being a subring of $M$.
\begin{equation}
\alpha(mm')\ =\ \alpha\1(m)\alpha\2(m')
\end{equation}
therefore multiplication $\mu_M\colon M\o_R M\to M$ is a left $A$-module map.
\begin{equation}
\alpha\triangleright 1_M\ =\ \eps_A(\alpha)1_M
\end{equation}
therefore the unit $_RR_R\to \,_RM_R$ is a left $A$-module map, too. This
means precisely that $M$ together with its ring structure, written as maps in
$_R\M_R$, is a monoid in $_A\M$, i.e., $M$ is a left $A$-module algebroid.

The invariants are determined as follows.
First of all $N \subset M^A$ is obvious. On the other hand if $m \in M^A$,
then $\beta_i(m) = \eps_A(\beta_i)m=\beta_i(1)m$ for each $i$, so for every
$\psi \in \mathcal{E} := \End M_N$,
$$ \psi(m) = \psi(b^1_i)b^2_i \beta_i(m)
= \psi(1)m, $$
whence also $\psi \circ \lambda(m')(m) = \psi(m'm) = \psi(m')m$
for each $m' \in M$. Thus, $\rho(m) \in \End {}_{\mathcal{E}}M = \rho(N)$,
and $m \in N$. Only here in the last step have we used that $M_N$ is
balanced.
\end{proof}
We set down some equivalent formulae for the invariant subring, the proof
of
which are left to the reader.
\begin{eqnarray}
M^A&:=&\{n\in M\,|\,\alpha\triangleright n=\eps_A(\alpha)n,\ \forall\alpha\in
A\,\}\\
&=&\{n\in M\,|\,\alpha\circ\rho(n)=\rho(n)\circ\alpha,\ \forall\alpha\in
A\,\}\\
&=&\{n\in M\,|\,\alpha\circ\lambda(n)=\lambda(n)\circ\alpha,\
\forall\alpha\in
A\,\}
\end{eqnarray}

\begin{exa}
\begin{rm}
If $M | N$ is an algebra extension with $N = k1$ trivial and $M$ finite
dimensional,
we recover Lu's bialgebroid $A = \End_k M$ \cite[3.4]{Lu}
since $R = M$. In case $R$ is not semisimple,
$A$ is a bialgebroid over $R$ which is not a weak bialgebra \cite{EN}.
This provides a wealth of examples of  action by bialgebroids.
\end{rm}
\end{exa}
\begin{rmk}
\begin{rm}
For each $\alpha \in A$, its coproduct $\cop_A(\alpha)$ may be considered a
map in $\Hom_{N-N}(M \o_N M, M)$ via
Prop. \ref{pro: tensor-square}.  We compute the simple form it takes:
\begin{equation}\label{eq: Lu}
\cop_A(\alpha)(m \o m') = \sum_i \gamma_i(m)c^1_i \alpha(c^2_im') =
\alpha(mm').
\end{equation}

\end{rm}
\end{rmk}

\begin{exa}
\begin{rm}
That $M_N$ should be balanced in the theorem is a necessary condition, for
consider $M$ to be the algebra of
$2$-by-$2$ matrices over a
field $k$ with $N$ the upper triangular matrices.  It is left as an exercise
to show that
$R$ is trivial ($k1_M$), $M |N$ is H-separable (therefore D2) since $1 \o 1
\in (M \o_N M)^M$, and
that $$\mathcal{E} := \End M_N  \cong M (\cong \End {}_NM \cong M \o_N M).$$
Consequently $A \cong R$ is trivial, so $M^A = M \neq N$.  But
$M_N$ is not balanced, since ${}_{\mathcal{E}}M = {}_MM$ and
$N \not \cong \End {}_MM = \rho(M)$.
It is also worth a mention that $M | N$ is not Frobenius, nor even
QF.
\end{rm}
\end{exa}

We next note that the endomorphism ring of a left D2 extension $M/N$
is isomorphic to a smash product of $M$ with the bialgebroid $A$.

\begin{cor}\label{cor: end = smash product}
$\End M_N$ is isomorphic to a smash product ring $ M \rtimes A$
via $m \rtimes \alpha \mapsto \lambda(m) \alpha$.
\end{cor}
\begin{proof}
Define $\pi: M \rtimes A \to \End M_N$ by the mapping just given.
By Proposition~\ref{pro: end}, $\pi$ is a linear isomorphism
of $M \o_R A$ with $\End M_N$.
  We compute using Eq.\ (\ref{eq: copA})  for $x \in M$:
\begin{eqnarray*}
\pi((m \rtimes \alpha)( m' \rtimes \beta))(x) & = & \pi(m \alpha\1(m') \rtimes
\alpha\2 \beta)(x) \\
& = & m \alpha\1(m')\alpha\2(\beta(x))  \\
& = & m \alpha(m'\beta(x)) \\
& = & \pi(m \rtimes \alpha) \circ  \pi(m' \rtimes \beta)(x)
\end{eqnarray*}
Hence, $\pi$ is a ring isomorphism.
\end{proof}

\begin{rmk}
\begin{rm}
Taking into account the corollary and the
finite projectiveness of $A$ over the  centralizer $R$ established in
the previous section, we propose that a D2 extension
 $M$ is an $B$-Galois extension of $N$
if $M$ is a balanced $N$-module.
 The justification for this terminology
 requires further investigation in a future paper.
\end{rm}
\end{rmk}


\section{The Right Bialgebroid $B$}\label{sec: B}

In this section, we define a right bialgebroid structure on $B=(M\o_NM)^N$
based on the restriction of the Sweedler $M$-coring on $M \o_N M$. However,
the proof will again depend on a left or right quasibasis.  In case
$N$ is trivial, the right
bialgebroid $B$ is Lu's left bialgebroid on the tensor-square up to a twist
by the antipode.
  The ring
stucture on $B$ is induced from an isomorphism $B \cong \End {}_M(M \o_N M)_M$.
There is a natural right action of $B$ on $\End {}_NM$ with fixed points
isomorphic to $M^{\rm op}$.

Let $B=(M\o_NM)^N$ the elements of which are denoted $b=b^1\o b^2$
suppressing
a possible summation.
$B$ is a ring with multiplication $bb'={b'}^1b^1\o b^2{b'}^2$ and unit $\one=1\o
1$. This multiplication does not extend to $M\o_NM$ but $M\o_NM$ is a left
$B$-module via
\begin{equation}\label{eq: kn}
b\cdot (m\o m')=mb^1\o b^2m'.
\end{equation}
The so defined ring
homomorphism $B\to\End_{M-M}(M\o_NM)$ is in fact an isomorphism. The inverse
is provided by $f\mapsto f(1\o 1)$.

Let $R$ be the centralizer of $N$ in $M$, $R=C_M(N)$. Define the ring
homomorphisms
\begin{eqnarray}
s_B\colon &R&\to \ B,\qquad s_B(r)=1\o r\,,\\
t_B\colon &R^{op}&\to\ B,\qquad t_B(r)=r\o 1\,.
\end{eqnarray}
Since we are going to make $B$ ino a right bialgebroid over $R$ we define
its $R$-$R$-bimodule via the actions
\begin{equation}
r\cdot b\cdot r'\ =\ bt_B(r)s_B(r')\ =\ rb^1\o b^2r'\,.
\end{equation}

\begin{lem}
Let $N\to M$ be a left D2 extension of rings.
Then the tensor product bimodule $B\o_R B$ is isomorphic, as an
$R$-$R$-bimodule, to $(M\o_N M\o_N M)^N$ where the bimodule structure of the
latter is defined by $r\cdot (m\o m'\o m'')\cdot r'=rm\o m'\o m''r'$. An
isomorphism is given by
\begin{equation}
\imath\colon B\o_R B\to (M\o_NM\o_N M)^N\,,\quad
b\o b'\mapsto b^1\o b^2{b'}^1\o {b'}^2\,.
\end{equation}
\end{lem}
\begin{proof}
That $\imath$ is a bimodule map is clear. To show that it is an isomorphism
we
write down its inverse using the left D2 quasibasis $\{b_i,\beta_i\}$ of
Lemma
\ref{lem: D2}.
\[
\imath^{-1}(t)=\sum_i b_i\o_R(\beta_i(t^1)t^2\o_N t^3)\,,\quad
t\in (M\o_NM\o_NM)^N \qed
\]
\renewcommand{\qed}{}\end{proof}

Now the right bialgebroid structure on the ring and bimodule $B$ is defined
by the following coproduct and counit
\begin{eqnarray}
\cop_B(b)&=&\sum_i (b_i^1\o_N b_i^2)\ \o_R\ (\beta_i(b^1)\o_N b^2)\\
\eps_B(b)&=&b^1b^2
\end{eqnarray}
By the lemma, $\cop_B(b) = \imath^{-1}(b^1 \o 1 \o b^2)$.

\begin{thm}
Let $N\to M$ be a D2 extension of rings. Then
$\bra B,R,s_B,t_B,\cop_B,\eps_B\ket$ is a right bialgebroid and $\End\,_NM$ is
a right $B$-module algebroid w.r.t. the action $\xi\triangleleft
b:=b^1\xi(b^2 -)$. The subring of invariants is $\rho(M)$, the right
multiplications with elements of $M$.
\end{thm}

\begin{proof} At first we check the bialgebroid axioms:

Coassociativity:
Apply $\imath_3\colon B\o_R B\o_R B\stackrel{~}{\to} (M\o_N M\o_N M\o_N
M)^N$,
$b\o b'\o b''\mapsto b^1\o b^2{b'}^1\o {b'}^2{b''}^1\o {b''}^2$ to both
hand sides of
$(\cop_B\o\id_B)\circ\cop_B=(\id_B\o\cop_B)\circ\cop_B$ and check that the
result on a $b\in B$ is $b^1\o 1\o 1\o b^2$ in both cases.  (The inverse
$\imath_3^{-1}$ sends $t \in  (M\o_N M\o_N M\o_N M)^N$ into
$$ \sum_{i,j} b_i \o_R (\beta_i(t^1)t^2 \o_N t^3\gamma_j(t^4)) \o_R c_j \in
B\o_R B\o_R B.)$$

Counit properties: Obvious.

The image of $\cop_B$ is in the Takeuchi $\x_R$ -product, $\cop_B(B)\subset
B\x_R B$:
\begin{eqnarray*}
(s_B(r)\o 1)\cop_B(b)&=&(b_i^1\o rb_i^2)\o_R(\beta_i(b^1)\o b^2)=\\
&=&\imath^{-1}\left(b_i^1\o rb_i^2\beta_i(b^1)\o b^2\right)
=\imath^{-1}(b^1\o r\o b^2)
\end{eqnarray*}
and similarly
\[
(1\o t_B(r))\cop_B(b)\ =\imath^{-1}(b^1\o r\o b^2)
\]

$\cop_B$ is multiplicative:
\begin{eqnarray*}
\cop_B(b)\cop_B(b')&=&(b_j^1b_i^1\o
b_i^2b_j^2)\o_R(\beta_j({b'}^1)\beta_i(b^1)\o
b^2{b'}^2)\\
&=&\imath^{-1}\left(b_j^1b_i^1\o b_i^2b_j^2\beta_j({b'}^1)\beta_i(b^1)\o
b^2{b'}^2\right)\\
&=&\imath^{-1}\left({b'}^1b_i^1\o b_i^2\beta_i(b^1)\o b^2{b'}^2\right)
=\imath^{-1}({b'}^1b^1\o 1\o b^2{b'}^2)\\
&=&\cop_B(bb')
\end{eqnarray*}

Unitalness: $\cop_B(\one)=\one\o\one$, $\eps_B(\one)=1$ are obvious.

$\eps_B$ is compatible with multiplication in the sense of axiom (v):
\[
\eps_B(t_B(\eps_B(b))b')=\eps_B((b^1b^2\o 1)({b'}^1\o
{b'}^2))={b'}^1b^1b^2{b'}^2=\eps_B(bb')
\]
and the same for $t_B$ replaced by $s_B$.

This finishes the proof that $B$ is a bialgebroid.

Module algebroid properties:
\[
(\xi\circ\xi')\triangleleft b=b_i^1\xi(b_i^2\beta_i(b^1)\xi'(b^2 -))=
(\xi\triangleleft b\1)\circ(\xi'\triangleleft b\2)
\]
The induced bimodule structure on $\End\,_NM$ is $r\cdot\xi\cdot
r'=\lambda(r)\circ\xi\circ\lambda(r')$
\[
\id_M\triangleleft b=\lambda(\eps_B(b))
\]

The invariants:
\[
(\End\,_NM)^B:=\{\xi\,|\,\xi\triangleleft b=\lambda(\eps_B(b))\circ\xi\,\}
\]
Clearly $\xi$ is an invariant iff
\[
b^1\o\xi(b^2m)=b_i^1\o b_i^2\beta_i(b^1)\xi(b^2m)=b_i^1\o
b_i^2\beta_i(b^1)b^2
\xi(m)=b^1\o b^2\xi(m)
\]
for all $m\in M$, $b\in B$. Thus
\[
1\o\xi(m)=\gamma_j(m)c_j^1\o \xi(c_j^2)=\gamma_j(m)c_j^1\o c_j^2\xi(1)
=1\o m\xi(1)
\]
Applying multiplication $\xi(m)=m\xi(1)$ follows. Thus an invariant $\xi =
\rho(\xi(1))$ and
belongs to  $\rho(M)$. The opposite inclusion is trivial. This proves
\[
(\End\,_NM)^B\ =\ \rho(M) \qed
\]
\end{proof}

Recalling the theory of the dual of a left bialgebroid in Section~2.6, we
have:

\begin{cor}\label{cor: core}
$B$ is isomorphic as bialgebroids over $R$ to the right bialgebroid dual
of $A$ via
the isomorphism in Eq.\ (\ref{B dual to A}). Similarly, $B$ is
isomorphic to the left bialgebroid dual ${}^*A$ via $\psi$ in
Theorem~\ref{thm: D2}.
\end{cor}
\begin{proof}
We prove the first statement and leave the second as an exercise.
Recall the nondegenerate pairing $\bra b, a \ket = b^1 a(b^2)$.
Let $A^*$ denote the right bialgebroid dual of $A$ with $\eta: B \to A^*$
the linear isomorphism given by $\eta(b) = \bra b, - \ket$.
We note that $\eta$ is an $R$-$R$-bimodule homomorphism,
since $$ \bra r \cdot b \cdot r', a \ket = rb^1 a(b^2 r') = r \bra b, a t(r')
\ket. $$
$\eta$
is a ring homomorphism since
$$
\bra bb',a \ket = {b'}^1 b^1 a(b^2 {b'}^2)
$$
while
$$
\bra b', \bra b, a\1 \ket \cdot a\2 \ket = {b'}^1 b^1 a\1(b^2) a\2({b'}^2) =
{b'}^1 b^1 \alpha(b^2 {b'}^2).
$$

$\eta$ is a homomorphism of corings since
\begin{eqnarray*}
\bra b, aa' \ket & = &  b^1 aa'(b^2) \\
& = & \sum_i b^1_i a( b^2_i \beta_i (b^1) a'(b^2)) \\
& = & \sum_i \bra b_i \beta_i(b^1) a'(b^2), a \ket = \bra b\1 \cdot \bra b\2,a'
\ket, a \ket. \qed
\end{eqnarray*}
\renewcommand{\qed}{}\end{proof}
\begin{rmk}
\begin{rm}
There is also a right action of $B$ on $\End M_N$ given by $\xi\triangleleft
b:=\xi(- b^1)b^2$. It however satisfies
\[
(\xi\circ\xi')\triangleleft b=(\xi\triangleleft b\2)\circ(\xi'\triangleleft
b\1)
\]
Its invariants are also the right multiplications with elements of $M$.
\end{rm}
\end{rmk}
\begin{rmk}
\begin{rm}
The coring $\bra B,R, \Delta_B, \eps \ket$ is a restriction of the Sweedler
coring \cite{Sw}
$$\bra M \o_N M,\, M ,\, \Delta:M \o_N M \to M \o_N M \o_N M ,\, \eps: M \o_N
M
\to M \ket.$$

If $N = k1$ for a ground field $k$ with $M$ finite dimensional, we recover
Lu's bialgebroid $B = M^{\rm op} \o_k M$ \cite[3.1]{Lu} up to a twist $S$.
In this
case, $B$ is a Hopf algebroid, with antipode $S$.
\end{rm}
\end{rmk}


\section{The Frobenius Case}
Depth two for Frobenius extensions is an important case to consider since
the Pimsner-Popa orthonormal
basis result shows that
finite index subfactors are Frobenius extensions \cite{GHJ,JS,EK,K}.
After recalling the basic construction or endomorphism ring theorem for
Frobenius extensions in Proposition~\ref{pro: basic construction}, we prove
that the most usual classical definition of depth two in terms of a Frobenius
extension $M | N$ is a special case of  depth two for ring extensions. We
show that $A = \End {}_NM_N$ and $B = (M \o_N M)^N$ are Frobenius extensions
over thecentralizer $R$. We introduce isomorphisms $\psi_A$, $\psi_B$
between $A$ and $B$ and the step two centralizers in the Jones
tower above $M | N$ denoted $\hat{A}$ and $\hat{B}$, respectively.
We prove an endomorphism ring theorem for D2 Frobenius extensions.

Recall that a ring extension $M |N$ is \textit{Frobenius} if
there is (a Frobenius homomorphism) $E \in \Hom_{N-N}(M,N)$
and (dual bases) $x_i, y_i \in M$ such that $\sum_i \lambda(x_i) E
\lambda(y_i)$
$= \id_M = \sum_i \rho(y_i) E \rho(x_i) $. Throughout
this section and part of the next, we assume $M | N$ is Frobenius with this
data.
We recall several facts about $M| N$.

\begin{pro}\label{pro: basic construction}
We have $\End M_N \cong M \o_N M$, which is a Frobenius extension over
$\lambda(M)\cong M$,
with Frobenius homomorphism $E_M = \mu$ and dual bases $\{ x_i \o 1 \}$, $\{
1 \o y_i \}$.
Moreover, $\End\,_NM$ and $\End M_N$ are anti-isomorphic.
\end{pro}
\begin{proof}
The isomorphism $\mathcal{F}: \End M_N \to M \o_N M$ is given by $f \mapsto
\sum_i f(x_i) \o y_i$
with inverse $m \o m' \mapsto \lambda(m) E \lambda(m')$.  A multiplication
on $M \o_N M$ is induced from composition of endomorphisms,
the $E$-multiplication given by $(m \o m')(m'' \o m''')
= m E(m'm'') \o m'''$ and unity $1_{1} = \sum_i x_i \o y_i$.
An anti-isomorphism is then given by
\begin{equation}
\label{phi}
\phi:\ \End\,_NM \stackrel{\cong}{\to} M \o_N M,\ \  f \mapsto \sum_i x_i \o
f(y_i)
\end{equation}
with inverse $m \o m' \mapsto \rho(m') E \rho(m)$.
The rest of the proof is
somewhat standard \cite{K}.
\end{proof}

We set $e_1 = 1 \o 1$ and  $M_1 := Me_1M = M \o_N M$.
Note that $M_1 \cong \End M_N$ via the $M$-$M$ map
induced by $e_1 \mapsto E$ and $M$ identified with $\lambda(M)$.  Note too
the key
identities
\begin{equation}
e_1 m e_1 = e_1 E(m) = E(m) e_1, \ \ \ E_M(me_1m') = mm'.
\end{equation}
In this notation $\{ x_i e_1 \}$, $\{ e_1 y_i \}$ are dual bases for $E_M \in
\Hom_{M-M}(M_1, M)$.

We note then that $A = \End\,_NM_N \cong (M \o_N M)^N$ via $\alpha \mapsto
\sum_i \alpha(x_i)
\o y_i$.

If we iterate this (basic) construction, we construct $M_2 = M_1 e_2 M_1$
with $e_2 m_1 e_2$ $= e_2 E_M(m_1) = E_M(m_1) e_2$ and $E_{M_1}(m_1^1 e_2
m_1^2) = m_1^1 m_1^2$ for each $m_1 \in M_1$.  Note that $M_2 = M_1 \o_M M_1
\cong M \o_N M \o_N M$.    We arrive at a generalized Jones tower,
$$ N \to M \into M_1 \into M_2 \into \cdots  $$
with \textit{Temperley-Lieb generators} $e_i \in M_i$ such that
\begin{equation}
e_i e_{i+1} e_i = e_i 1_{M_{i+1}}, \ \ e_{i+1} e_i e_{i+1} = e_{i+1}, \ \ e_i
e_j = e_j e_i
\end{equation}
if $| i - j | > 1$. Note that the $e_i$ are not the Jones projections
even if they exist, $e_i^2 \neq e_i$. For example if $M|N$ is not a split
extension then there is no unit preserving Frobenius homomorphism $E$.
However, the Temperley-Lieb generators exist for any Frobenius extension as
shown above.

We also have the Pimsner-Popa relations:
\begin{equation}
m_i e_i = E_{M_{i-1}}(m_ie_i)e_i, \ \ e_im_i = e_i E_{M_{i-1}}(e_im_i)
\end{equation}
for $m_i \in M_i$, $i = 1,2$ and $M_0 := M$.

Introduce the notation $\hat{A} := M_1^N$ for the centralizer $C_{M_1}(N)$
in $N \to M \into M_1$.  We introduce the canonical isomorphism $\psi_A$
of $A = \End
{}_NM_N$  with $\hat{A}$ given by the restriction of $\mathcal{F}$ above to $A$:
$$\alpha \longmapsto \sum_i \alpha(x_i) e_1 y_i, $$
with inverse $a^1 e_1 a^2 \mapsto \lambda(a^1) \circ E \circ \lambda(a^2). $
Similarly, let $\hat{C} := M_2^N$, which is isomorphic as rings to the step three
centralizer $C$ introduced in Section~3.

We first show that classical depth two extensions are depth two in the sense
of this paper.  It is known that a semisimple pair $N \subset M$
over a field, where $M$ has faithful
trace $T$ that restricts to a faithful trace on $N$, is a (split,
separable) Frobenius
extension; cf.\ \cite[Prop.\ 2.6.2]{GHJ}. Also subfactors of finite
index are Frobenius extensions by the Pimsner-Popa orthonormal basis
result \cite{GHJ, K}; these have semisimple centralizers.
Of course, a module
over a semisimple ring is always projective.

\begin{pro}\label{pro-classical D2}
Suppose $M | N$ is Frobenius extension, $\hat{A}_R$ and ${}_R\hat{A}$
are f.g. projective,  and $\hat{C} = \hat{A}e_2 \hat{A}$.
Then $M | N$ is a depth two ring extension.
\end{pro}
\begin{proof}
We first show that $E_M: M_1 \to M$ has dual bases in $\hat{A}$.
 By the classical D2 hypothesis on $\hat{C}$,
$1_{M_2} =  \sum_k a_k e_2 b_k$ for some $a_k, b_k \in \hat{A}$.
  Let $m_1 \in M_1$, then:
$$ e_2 m_1 = \sum_k e_2 m_1 a_k e_2 b_k = \sum_k e_2 E_M(m_1 a_k) b_k. $$
By applying $E_{M_1}$ we arrive at $m_1 = \sum_k E_M(m_1 a_k)b_k$;
similarly,
\newline
 $m_1 = \sum_k a_k E_M(b_k m_1)$, so $\{ a_k \}$, $\{ b_k \}$
are indeed dual bases for $E_M$.

It follows that $M_1 \cong M \o_R \hat{A}$ as $M$-$N$-bimodules via
$m_1 \mapsto \sum_k E_M(m_1 a_k) \o b_k$ with inverse
mapping given simply by $m \o a \mapsto ma$: note that $E_M(\hat{A}) \subseteq
R$. Since $M \o_N M \cong M_1$ and ${}_R\hat{A}$ is f.g. projective,
it follows that $M | N$ is right D2.  Similarly,
$M_1 \cong \hat{A} \o_R M$ as $N$-$M$-bimodules and
$M | N$ is left D2.
\end{proof}

\begin{pro}\label{pro: knD2}
If $M|N$ is a left or right D2 Frobenius extension, then $E_M: M_1 \to M$ has
dual bases in $\hat{A}$.
\end{pro}
\begin{proof}
Let $b_i \in B, \beta_i \in A$ be a left D2 quasibasis, then
$\{ b_i \}$, $\{ \sum_j \beta_i(x_j) e_1 y_j \}$ are dual bases, obviously in
$M_1^N$,
for $E_M$. As a matter of fact
\begin{eqnarray*}
\sum_{i,j} E_M(me_1m' b_i^1 e_1 b^2_i) \beta_i(x_j) e_1 y_j &=&
\sum_{i,j} mE(m'b^1_i)b^2_i \beta_i(x_j) e_1 y_j \\
&=& m E(m'x_j) e_1 y_j = m e_1 m'
\end{eqnarray*}
for $m, m' \in M$, and
$$ \sum_{i,j} b^1_i e_1 b^2_i E_M(\beta_i(x_j) e_1 y_j m e_1 m') = \sum_i
b_i^1 e_1 b_i^2
\beta_i(m)m' = m e_1 m'. $$
The proof starting with a right D2 quasibasis is similar.
\end{proof}

\begin{cor}
We have $M_1 \cong M \o_R \hat{A}$ via $m \o a \mapsto ma$ for each $m \in M,
a \in \hat{A}$.
\end{cor}
\begin{proof}
An inverse is given by $m_1 \mapsto \sum_{i,j} E_M(m_1 b_i) \o_R
\beta_i(x_j)e_1 y_j
\in M \o \hat{A}$ by Proposition~\ref{pro: end}.
\end{proof}

Similarly we show $M_1 \cong \hat{A} \o_R M$ via $a \o m \mapsto am$.
If $R$ is a field coincident with centralizers of $M$ and $N$
as in \cite{KN}, it follows from the proposition that $M_1$ is f.g. free
as a left or right natural $M$-module.

\begin{cor}\label{cor: smash product}
$M_1$ is isomorphic to a smash product algebra:  $M_1 \cong M \rtimes A$.
\end{cor}
\begin{proof}
Define $\Pi: M \rtimes A \to M_1$ by $\Pi(m \rtimes \alpha) = \sum_i m
\alpha(x_i)e_1 y_i$
for all $\alpha \in A, m \in M$.  We see then that $\Pi$ is a
composition of two algebra isomorphisms,
$\pi$ in Corollary~\ref{cor: end = smash product}
and $\mathcal{F}$ above (cf.\ Proposition~\ref{pro: basic construction}).
\end{proof}

 From Section~3 we recall the step~3 centralizer
$C = \End_{N-M} (M \o_N M)$.

\begin{cor}\label{cor: frobenius}
If $M | N$ is D2 Frobenius, then the ring extensions
$R \into A$, $r \mapsto \lambda(r)$
and $C | A$ (given by Eq.\ (\ref{eq: nn})) are Frobenius extensions.
\end{cor}
\begin{proof}
If $M | N$ is Frobenius, we arrive at $A | R$ Frobenius
 from the proposition by restriction of $E_M$ to $A$ (identified
with $\{ \alpha(x_i)e_1y_i|\, \alpha \in A \}$)
and noting that $E_M(A) \subseteq R$.  Since $C \cong \End {}_R A$, we
conclude
 from the (left) endomorphism ring theorem for Frobenius extensions
\cite{K}
that $C | A$ via right regular representation is Frobenius.
\end{proof}

Conversely, $A | R$ Frobenius implies
$\mathcal{E}| M$  is Frobenius by Prop.\ \ref{pro: end}.  If $M_N$ is
a progenerator, then a endomorphism ring theorem-and-converse
assures us that $M | N$ is also Frobenius (cf.\ \cite{K}). By the same
token, $A | R$ is Frobenius iff $C | A$ since ${}_RA$ is a progenerator
(Theorem~\ref{thm: D2}).

The next result is an endomorphism ring theorem for Frobenius D2 extensions,
and answers a question posed at the end of \cite{KN}.

\begin{thm}
\label{thm-endo}
If $M | N$ is a left D2
extension, then $M_1 | M$ is a right D2  extension.
Similarly,
if $M | N$ is right D2, then $ M_1|M$ is left D2.
\end{thm}
\begin{proof}
We note the bimodule ${}_{M_1}M_N$ given by
\begin{equation}
(m e_1 m') \cdot m'' \cdot n = m E(m'm'')n,
\end{equation}
which is of course isomorphic to the natural bimodule ${}_{\mathcal E}M_N$
where
$\mathcal{E} = \End M_N$. Now tensor from the left the first isomorphism in
Def.~\ref{defi: D2}
by this bimodule:

$$ {}_{M_1}M \o_N M \o_N M_M \oplus * \cong \oplus^n {}_{M_1}M \o_N M_M  $$
which is isomorphic to
$$ {}_{M_1}M_1 \o_M M_1{}_M \oplus * \cong \oplus^n {}_{M_1}M_1{}_M, $$
the condition for $M_1 | M$ to be right D2. The second statement is proven
similarly.
\end{proof}

Let $\hat{B} := M_2^M$, and note the canonical algebra isomorphism $\psi_B: B
\cong \hat{B}$
given by
$$\psi_B( b) = \sum_i x_i b^1 e_1 b^2 e_2 e_1 y_i $$
with inverse $$b^1 e_2 b^2 \mapsto b^1 E_M(b^2 e_1) $$
($b^1, b^2 \in M_1$)
obtained by following the ring isomorphisms
\begin{equation}\label{eq: bee}
M_2^M \cong \Hom_{M-M}(M_1,M_1) \cong (M \o_N M)^N
\end{equation}
via first the Frobenius map $\Psi: x \o y \mapsto \lambda(x) E_M \lambda(y)$,
for each $x,y \in M_1$,
composed with the general map $\Phi: f \mapsto f(1_M \o 1_M)$ for each $f \in
\End\,_M(M_1)_M$.
The $E_M$-multiplication on $M_2^M$ is therefore identifiable with
composition as well as
the multiplication on $B$ from Section~5.

Now the endomorphism ring theorem for D2 Frobenius
extensions may be used to show, in a
similar way to the earlier propositions in this section, that $E_{M_1}$
has dual bases in $\hat{B}$,
$B | R^{\rm op}$ is a Frobenius extension
and $M_2$ is  a smash product of $M_1$ and $\hat{B}$.


\section{The Biseparable Case}

Suppose $M | N$ is a Frobenius D2 $R$-algebra extension
where the centralizer $R$ is \textit{trivial},
i.e., coincides with  the centers of $M$ and $N$.
In this case,
$A$ and $B$ are bialgebras which are finitely generated projective
over $R$ by Cor.\ \ref{cor: frobenius} and its analog for $B$.
In this section, under
the additional constraints that $R$ coincides with
a ground field and $M | N$ is a biseparable algebra extension,
we show (by means of anti-isomorphisms $\phi_A$ and
$\phi_B$) that $A$ and $B$ are isomorphic as bialgebras with
the
bialgebra structures on $\hat{A}$ and $\hat{B}$ defined
in \cite{KN}.  As a consequence,
$A$ and $B$ are dual semisimple Hopf algebras
 with
Galois (Ocneanu-Szyma\'{n}ski) actions
 on $M$ and $M_1^{\rm op}$, respectively.  Finally,
we drop the Frobenius assumption on $M | N$ and show by means
of the techniques in Section~3 that
a D2 biseparable extension is quasi-Frobenius (QF).

The next proposition shows that $A$ acts on $M$ via the same formula as the
Hopf algebra action of $\End_{M-M} M_1$ on $M_1$ in \cite[Eq. (21)]{KN}. The
proof
does not make use of the triviality assumption on $R$.

\begin{pro}\label{pro: Ocneanu-Szymanski action}
If $M | N$ is Frobenius, then the action of $A$ on $M$ defined in
Section~4 is
the Ocneanu-Szyma\'nski action,
$$a \triangleright m = E_M(am e_1).$$ The algebra homomorphism $\pi: M \rtimes
A
\to \End M_N$ given by $m \rtimes a \mapsto \lambda(m) (a \triangleright -)$
is an isomorphism which
fits into a commutative triangle with the isomorphisms given in
Cor.\ \ref{cor: smash product} and Prop.\ \ref{pro: basic construction}.
\end{pro}

\[
\parbox[c]{4in}{
\begin{picture}(180,70)
\put(20,50){$M\rtimes A$}
\put(55,52){\vector(1,0){85}}  \put(95,57){$\pi$}  \put(95,42){$\cong$}
\put(150,50){$\End M_N$}
\put(37,42){\vector(2,-1){50}} \put(50,21){$\Pi$}
\put(92,10){$M_1$}
\put(160,42){\vector(-2,-1){50}} \put(138,21){$\mathcal{F}$}
\end{picture}
}
\]

\begin{proof}
We let $a \in (M \o_N M)^N$ be the image of $\alpha \in A$ under the
isomorphism
above.  Then:
$$ E_M(am e_1) = E_M(\alpha(x_i)e_1y_ime_1) = E_M(\alpha(x_i) E(y_i m) e_1) =
\alpha(m).  $$

The commutativity of the diagram is immediate from the definitions.
\end{proof}
\vskip 0.5truecm
The bimodule action on $A$ induced by $R \subset M \into M_1$ is given by the
somewhat different
formula $r \cdot \alpha \cdot r' = \lambda(r) \alpha \lambda(r')$ (cf.\
Eq.\ (\ref{eq: alf})).
However, the two bimodule structures coincide when $R$ is trivial.

We introduce two canonical anti-isomorphisms  of $ \phi_A: A \to \hat{A} =
M_1^N$
and $\phi_B: B \to \hat{B} = M_2^M$ given by
$$ \phi_A(\alpha) = \sum_i x_i e_1 \alpha(y_i), \ \
\phi_B(b) = \sum_i x_i e_1 e_2 b^1e_1b^2 y_i.
$$

The Ocneanu-Szyma\'nski action $\triangleright'$  of $\hat{B}$ on $M_1$ in
\cite[Eq. (21)]{KN} is related to our action of $B$
on $\mathcal{E}' = \End\,_NM$ via the anti-isomorphism $\phi$ in
Eq.\ (\ref{phi}) as we see next.  Again, we do not need triviality of $R$.

\begin{pro}
If $M | N$ is a depth two Frobenius extension, then for
every $b \in B$, $f \in \mathcal{E}'$:
$$\phi_B(b) \triangleright' \phi(f) :=
E_{M_1}(\phi_B(b) \phi(f) e_2) = \phi(f \triangleleft b). $$
\end{pro}
\begin{proof}
\begin{eqnarray*}
E_{M_1}(\phi_B(b) \phi(f) e_2) & = &
\sum_{i,j} E_{M_1}( x_i  e_1  e_2  b^1 e_1 b^2 y_i x_j e_1 f(y_j) e_2) \\
           & =& \sum_i E_{M_1}(x_i  e_1  E_M(b^1e_1 E(b^2 y_i
x_j)f(y_j))e_2) \\
          & = & \sum_i x_i e_1 b^1 f(b^2 y_i) = \phi( f \triangleleft b).
\qed
\end{eqnarray*}
\renewcommand{\qed}{}\end{proof}

A ring extension $M | N$ is said to be \textit{biseparable} if $M_N$ and
${}_NM$ are f.g. projective while $M | N$ is a
separable extension (i.e.,
$\mu: M \o_N M \to M$ is split $M$-$M$-epimorphism) and a split extension
(i.e., there is $N$-bimodule $V$ such that $M \cong N \oplus V$ as
$N$-bimodules).  If $R$ is trivial, a biseparable Frobenius extension
of $R$-algebras coincides with the notion of strongly separable extension
\cite{K,KN}.

\begin{thm}\label{thm: HA iso}
If $R$ is a field and $M | N$
is a biseparable Frobenius $R$-algebra extension of depth two, then
$A$ and $B$ are dual semisimple Hopf algebras isomorphic to
$\hat{A}$ and $\hat{B}$, respectively.
\end{thm}
\begin{proof}
Since $R$ is trivial, $A$ and $B$ are dual bialgebroids over $R$,
it follows easily that $A$ and $B$ are dual bialgebras.
We next note that the nondegenerate
pairing $\bra b,a \ket = b^1 a(b^2)$ is equal to the nondegenerate pairing
$$\bra \phi_A(a), \psi_B(b) \ket' :=  E_M E_{M_1}(\psi_B(b)e_1 e_2
\phi_A(a))$$
analyzed in \cite[4.4]{KN}, since:
\begin{eqnarray*}
E_M E_{M_1}(\psi_B(b)e_1 e_2 \phi_A(a)) & = &   \sum_{i,j}
E_M E_{M_1}(x_i b^1 e_1 b^2 e_2 e_1 y_i e_1 e_2 x_j e_1 a(y_j)) \\
& = & \sum_{i,j}  E_M E_{M_1}(x_i b^1 e_1 b^2 e_2 E_M(e_1 E(y_i))x_je_1
a(y_j)) )\\
& = & \sum_{i,j}  E_M(x_i b^1 e_1 b^2 E(y_i) x_j e_1 a(y_j))\\
& = & \sum_j b^1 E(b^2 x_j) a(y_j) = \bra b,a \ket.
\end{eqnarray*}

In \cite[Section 4]{KN}, it is shown that for $a' \in \hat{A}, b' \in
\hat{B}$
$$ \bra a', b' \ket' = E_M E_{M_1} ( a' e_2 e_1 S(b')) := \bra a',S(b')
\ket'' $$
for antipode $S: \hat{B} \to \hat{B}$.
Now let $\psi_B(b) = b'$ and $\alpha, \alpha' \in A$.  We compute that
$\psi_B$ is a coalgebra homomorphism:
\begin{eqnarray*}
\bra \phi_A(\alpha'),\psi_B(b\1) \ket'
\bra \phi_A(\alpha), \psi_B(b\2) \ket' & = & \bra b\1, \alpha' \ket \bra b\2,
\alpha \ket \\
& = & \bra b, \alpha'\alpha \ket \\
& = & \bra \phi_A(\alpha) \phi_A(\alpha'), S(b') \ket'' \\
& = & \bra \phi_A(\alpha), S({b'}\2) \ket'' \bra \phi_A(\alpha'), S({b'}\1)
\ket'' \\
& = &  \bra \phi_A(\alpha'), {b'}\1 \ket' \bra \phi_A(\alpha),{b'}\2 \ket',
\end{eqnarray*}
since $S$ is coalgebra anti-isomorphism and by definition of $\Delta(b')$
in \cite{KN}. Finally, $\psi_B$ preserves the counit:
\begin{eqnarray*}
\eps_{\hat{B}}(\psi_B(b)) &=& \sum_i E_{M_1}(e_2x_i b^1 e_1 b^2 e_2 e_1 y_i)
\\
& = & \sum_i E_{M_1}(e_2 E_M(x_i b^1 e_1 b^2) e_1 y_i) \\
& = & \sum_i x_i b^1 b^2 e_1 y_i = \eps_B(b)1_{M_1}
\end{eqnarray*}
by triviality of $R$ and \cite[3.13, 4.3]{KN}.
It follows that $\psi_B$ is a bialgebra isomorphism, whence
$B$ has an antipode. Semisimplicity of $A$ and $B$ follow
 from \cite{KN}.
$A$ is then also a Hopf algebra since it is the dual of $B$.
\end{proof}

Moreover, the antipode is involutive, $S^2 = \id$, by a powerful theorem
of Etingof and Gelaki \cite{EG}.
\subsection{D2 biseparable extensions are QF}

In this subsection, we no longer assume $M | N$ is Frobenius:
in fact, we will be interested in  \textit{when} D2 biseparable
extensions are Frobenius.
A \textit{depth one} ring extension $M | N$ is  a centrally projective ring
extension defined
in  Example~\ref{exa-cpre}: compare \cite[3.1]{KN}.  It follows from
Example~\ref{exa-cpre}
that a depth one extension is automatically D2.

The following theorem answers \cite[Problem 3.8]{CK} for
depth two extensions.
A ring extension $M | N$ is left QF (quasi-Frobenius)
if $M_N$ and ${}_NM$ are f.g. projective and ${}_NM^*_M
\oplus * \cong \oplus^n {}_NM_M$ \cite{Mu}.  Similarly there
is a notion of right QF extension with two-sided QF extensions
being denoted  by ``QF.'' Of course, a QF extension is a weakening
of the notion of Frobenius extension.

It is already well-known and easily derived
that a depth one (bi)separable extension is
QF (e.g., see \cite{CK}: in fact, it
is Frobenius \cite{Su}). The same is true of depth two extensions:

\begin{thm}
\label{thm-biseparable}
A depth two biseparable extension is QF.
\end{thm}
\begin{proof}
Since $M | N$ is separable, it follows easily that
$\mathcal{E} | M$ (identified with $\lambda(M)$)
is split with bimodule projection given by $f \mapsto \sum_i f(x_i)y_i$
where $\sum_i x_i \otimes y_i$ is a separability element.
Let ${}_MW_M$ be the complementary bimodule satisfying $\mathcal{E} \cong
M \oplus W$ as $M$-bimodules. By Corollary~\ref{cor-endo}
for left D2 extensions
we note that ${}_M\mathcal{E}_N \oplus * \cong \oplus^n {}_MM_N$.
Then also ${}_MW_N \oplus * \cong \oplus^n {}_MM_N$.
Since the canonical map $W \o_N M \to W$ is a split right $M$-epimorphism
by separability, it follows from ${}_NW \o_N M_M \oplus * \cong \oplus^n
{}_N M \o_N M_M$ and the definition of left D2 extension that
$$ {}_NW_M \oplus * \cong \oplus^m {}_NM_M.$$

On the other hand,
since $M | N$ is split, we have
$M \cong N \oplus V$ for $N$-bimodule
$V$, so $\mathcal{E} \cong \Hom(M_N, N_N \oplus V_N)$;
whence the isomorphism of $N$-$M$-bimodules,
$$ M \oplus W \cong \Hom(M_N,N_N) \oplus \Hom(M_N,V_N). $$

 From the two displayed equations, we conclude that $M^* \oplus * \cong
\oplus^{m+1} M$ as $N$-$M$-bimodules.  Hence $M | N$ is left
QF.  We prove similarly that a right D2 biseparable extension
is right QF.
\end{proof}

By comparing with the definition of depth three in \cite[3.1]{KN},
we propose
that
a ring extension $M | N$  be  \textit{right depth three} if
${}_{\mathcal{E}}M \o_N
M \o_N M_N$ and ${}_{\mathcal{E}}M \o_N M_N$ are H-equivalent modules;
and \textit{left depth three} if ${}_NM \o_N M \o_N M_{\mathcal{E'}}$
and ${}_NM \o_N M_{\mathcal{E'}}$ are H-equivalent.
We recall our notations $\mathcal{E} = \End M_N$ and $\mathcal{E'} =
\End\,_NM$.
The following is proved in the same way as Theorem~\ref{thm-endo}.

\begin{pro}
A depth two extension is depth three.
\end{pro}

We propose the following problem in extension of Theorem
\ref{thm-biseparable}:
is a biseparable depth three extension
QF or  Frobenius?   Yet another problem is to determine a reasonable
definition of
finite depth ring extensions.


\section{The Irreducible Case}

 From \cite{KN} we recall (and slightly extend the notion) that a
$K$-algebra
extension $N \to M$ is \textit{irreducible} if the centralizer
is trivial: $R = K1_M$ for $K$ a commutative ground ring.
In this section, we show that a depth two irreducible Frobenius extension
$M/N$
has \textit{Hopf algebras} $A$ and $B$ with bialgebra
structure defined as before. Because
of the results in Sections~6 and ~7, this extends by entirely different means
the main theorem 1.1 in \cite{KN}.  To be precise,
 we obtain the same results without the hypotheses of
 biseparability
and \textit{field} $K$; however, we introduce the new condition of
$M_N$ being balanced (as noted in Section~4)
to obtain the fact that $M^A = N$. Finally, we prove that
$B$, resp.\ $A$, is $K$-separable if and only if $M |N$
 is a separable, resp.\ split,  extension. For $K$ a characteristic
zero field, this implies that an irreducible D2 Frobenius extension
is split iff it is separable.

\begin{thm}
Suppose $M/N$ is a depth two irreducible Frobenius extension.
Then $A := \End {}_NM_N$ and $B := (M \o_N M)^N$ are dual Hopf
algebras acting on $M$ and $\End {}_NM$ respectively, with $M_1 \cong
M \rtimes A$.
\end{thm}

\begin{proof}
 From Section~4 recall that $A$ is a left bialgebroid over $R = K$;
whence a $K$-bialgebra which by Theorem~\ref{thm: D2} is a
progenerator $K$-module. We also recall that $A$ acts on $M$ with
$M_1$ isomorphic as rings to $M \rtimes A$.
 From Proposition~\ref{pro: rightdual} and
Corollary~\ref{cor: core},  $B$ is the bialgebra dual of $A$, since
$$ \bra bb',a \ket = \bra b, a\1 \ket \bra b',a\2 \ket, \ \
\bra b, aa' \ket = \bra b\1, a \ket \bra b\2, a' \ket.  $$
It suffices then to show that $A$ has an antipode.

Let $E: M \to N$ be a Frobenius homomorphism with dual bases $\{ x_i \}$,
$\{ y_i \}$, and $b_i, \beta_i$ be a left D2 quasibasis for $M/N$.
We now claim that $\psi: A \to K$ defined by $\psi(\alpha) = \sum_j
\alpha(x_j)y_j$is a Frobenius homomorphism satisfying
\begin{equation}
a\1 \psi(a\2) = \psi(a) 1_A
\end{equation}
for every $a \in A$.  $\psi$ is shown to be a Frobenius homomorphism  by
either noting that it corresponds to $E_M$ restricted to $\hat{A}$ via
the isomorphism $A \cong \hat{A}$ given in Corollary~\ref{cor: frobenius},
or computing that $\{ b_i^1E(b^2_i -) \}$, $\{ \beta_i \}$, are dual
bases for $\psi$.Now we compute:
\begin{eqnarray*}
 a\1 \psi(a\2) & = & a(-b^1_i)b_i^2 \beta_i(x_j)y_j \\
               & = & \alpha(-x_j)y_j \\
               &  = & \alpha(x_j)y_j \id_M = \psi(a) 1_A,
\end{eqnarray*}
since $\sum_i x_i \o y_i \in (M \o_N M)^M$, so $a(mx_j)y_j = \psi(a)m$
for $m \in M$.

We note next that $E$ is left norm for the augmented Frobenius algebra
$(A,\psi, \eps)$, since for each $a \in A$:
$$ \psi(aE) = \sum_j a(E(x_j))y_j = a(1) \sum_j E(x_j) y_j = \eps(a). $$
Now it follows from Eq.\ (\ref{eq: copA}) and a standard lemma
(due to Pareigis) that
\begin{equation}
\label{antipode formula}
S: A \to A, \ \ S(a) = E\1 \psi(a E\2) = \sum_{i,j} E(-b^1_i)b^2_i a
\beta_i(x_j)y_j
\end{equation}
is an antipode for $A$, since $\mathcal{A} := \Hom_K(A,A)$ is f.g.
projective algebra
with respect to the convolution product $*$ induced from $A$,
clearly $1_A * S = 1_{\mathcal{A}}$, whence $S * 1_A = 1_{\mathcal{A}}$.
\end{proof}

The theorem provides the key to computing formulas for the Hopf algebra
structureson $\hat{A}$ and $\hat{B}$ in \cite[Section~6]{KN}
and its extension to the nonbiseparable case
with commutative ground ring.  For example, we show that the action of
$\hat{A}$
on $M$ in \cite{KN}, expressed by a conjugation formula (Eq.\ 25),
 is indeed given by the Ocneanu-Szyma\'nski action. Let $\hat{S}$ be
the antipode induced from $S$ on $\hat{A}$ by the
Hopf algebra isomorphism $\psi_A(\alpha) = \sum_j \alpha(x_j) e_1 y_j$
(cf.\ Theorem~\ref{thm: HA iso}).

\begin{pro}
For $a \in \hat{A}$ and $m \in M$, we have
\begin{equation}
E_M(ame_1) = a\1 m \hat{S}(a\2)
\end{equation}
\end{pro}
\begin{proof}
Let $\psi_A(\alpha) = a$ for $\alpha \in A$.
By Proposition~\ref{pro: Ocneanu-Szymanski action}, it will suffice to compute
that $a\1 m \hat{S}(a\2) = \alpha(m)1_{M_1}$. We compute using Eq.\
(\ref{antipode formula}):
\begin{eqnarray*}
a\1 m \hat{S}(a\2) & = & \sum_{i,j,k,r,s}
\alpha(x_jb^1_i)b^2_i e_1 y_j m E(x_k b^1_r)b^2_r \beta_i
\beta_r(x_s)y_se_1 y_k \\
& = & \sum_{i,j,k,r,s} \alpha(x_jb^1_i)b^2_i e_1 E(y_j m E(x_k b^1_r)b^2_r
\beta_i \beta_r(x_s)y_s) y_k \\
& = & \sum_{k,r,i,s}\alpha( m E(x_k b^1_r)b^2_r b^1_i)b^2_i\beta_i
      \beta_r(x_s)y_s e_1 y_k \\
& = & \sum_{k,r,s} \alpha(m  E(x_k b^1_r)b^2_r\beta_r(x_s))y_s e_1 y_k \\
& = & \sum_{k,s} \alpha(m E(x_k x_s)) y_s e_1 y_k \\
& = & \alpha(m) \sum_k x_k e_1 y_k = \alpha(m) 1_{M_1},
\end{eqnarray*}
since $\sum_s \beta_i \beta_r(x_s)y_s \in R = K$.
\end{proof}

We next note a criterion for when an irreducible D2 Frobenius and
proper
extension $M/N$ is split (i.e., $N \oplus * \cong M$ as $N$-bimodules).

\begin{pro}\label{pro: split}
$M/N$ is split $\Leftrightarrow$ $A$ is $K$-separable.
\end{pro}
\begin{proof}
Since $M/N$ is Frobenius with Frobenius homomorphism $E: M \to N$, it is split
iff there is $d \in R$ such that $E(d) = 1$. Since $R$ is trivial, $d  \in K$,
so $E(1)d = 1$ and $ \eps(E) = E(1)$ is invertible.  Since $E$ is a left norm
in $A$, it is a left
integral, or by direct computation for $m \in M$ and $\alpha \in A$:
$$ \alpha E(m) = \alpha(1) E(m) = \eps(\alpha) E(m). $$
But then  $A$ is $K$-separable iff $\eps(E)$ is invertible
(e.g., \cite[5.2]{KS}).
\end{proof}

Similarly $M/N$ is separable  iff there is $d \in R$ such that
$\sum_i x_idy_i = 1$.
Then $E_M(1_{M_1}) = \sum_i x_iy_i$ is invertible in $K$.
Now the multiplication in $B$ yields
\begin{equation}\label{eq: right integral}
 1_{M_1} b = b^1 x_j e_1 y_j b^2 = 1_{M_1} b^1 b^2 = 1_{M_1} \eps(b),
\end{equation}
whence $1_{M_1}$ is a right integral for the Hopf algebra $B$.
We then similarly complete the proof of the next proposition:
 \begin{pro}\label{pro: separable}
$M/N$ is separable $\Leftrightarrow$ $B$ is $K$-separable.
\end{pro}

Next we conclude from Proposition~\ref{pro: Ocneanu-Szymanski action}
and \cite[1.1]{U} (cf.\ \cite[Fig.\ 2]{KN})
 that $M/N$ is a $B$-Galois extension of $K$-algebras,
a generalization of \cite[Theorem 6.5]{KN}.  We assume
that $M_N$ is balanced as used in Section~4.

\begin{cor}
An irreducible Frobenius  extension of depth two is a
Hopf-Galois extension.
\end{cor}

If $K$ is a field of characteristic zero, it follows from the Larson-Radford
theorem \cite{LR} that $A$ is a semisimple Hopf algebra
over $K$ iff its dual $B$ is semisimple.  From
the propositions above we deduce:
\begin{cor}
Suppose $K$ is a field of characteristic zero and $M/N$
is an irreducible D2 Frobenius $K$-algebra extension.  Then
$M$ is a split extension of $N$ if and only if $M$ is a separable
extension of $N$.
\end{cor}

\section{When $A$ and $B$ are Weak Hopf Algebras}

In this section we assume that $M | N$ is a depth two Frobenius extension
of algebras over a \textit{field} $K$
\textit{where the centralizer $R$ is a  separable $K$-algebra}.
We provide a new characterization of  separable algebra $R$
as an index one Frobenius algebra; i.e., possessing a Frobenius
system $\phi, e_i, f_i$ such that $\sum_i e_i f_i = 1$ and
$\sum_i \phi(re_i)f_i = r = \sum_i e_i \phi(f_i r)$ for each $r \in R$.
We prove that under this
condition on $R$  the step two centralizers $A$ and $B$ have dual weak Hopf
algebra structures.  (For weak Hopf algebra theory,
see \cite{BSz,BNS,KN2,N,NV3,Sz,V}.) It then follows from
Sections 2-6 that the action of $A$ on $M$ is the usual action
of a weak Hopf algebra as is  $\End M_N$ an ordinary
smash product of weak Hopf algebra with its module algebra \cite{NV3}.
This generalizes the results in \cite{KN2} by removing the Markov trace
on $M$, the conditions of biseparability and
 \textit{symmetric} Frobenius on $M | N$,
and relaxing the condition that $R$ be strongly separable.  However,
we note again that our approach yields
 $N = M^A$ by requiring that $M_N$ is a balanced module.
Finally we extend
Propositions~\ref{pro: split} and~\ref{pro: separable} to weak
Hopf algebras in
the case of ground fields.

First we begin with two useful lemmas for Frobenius algebras.

\begin{lem}
\label{lem: R-dual = K-dual}
Suppose $R$ is a Frobenius $K$-algebra and $V$ is a right $R$-module.
Then $\Hom_R(V,R) \cong \Hom_K(V,K)$.
\end{lem}
\begin{proof}
Let $\phi: R \to K$ be a Frobenius homomorphism with dual bases $\{ e_j \}$
and $\{ f_j \}$.  Then the mapping $\Hom_R(V,R) \to \Hom_K(V,K)$ given
by $f \mapsto   \phi \circ f$ has inverse given by $g \mapsto \sum_i
g(-e_i)f_i$.
\end{proof}

Recall that a linear functional $\phi: R \to K$ is \textit{left (right)
nondegenerate}
if $\phi(xR) = 0$ ($\phi(Rx) = 0$) implies $x = 0$.

\begin{lem}
\label{lem: left iff right nondegen}
If $K$ is a field and $R$ is a finite dimensional
 $K$-algebra, then $\phi: R \to K$
is left nondegenerate iff it is right nondegenerate.
\end{lem}
\begin{proof}
If $\phi$ is left nondegenerate, then by dimension comparison, we see
 $x \mapsto \phi x$ is an isomorphism of $A$ with its $K$-dual $A^*$.
Then we can find dual bases $\{ e_i \}$, $\{ \phi f_i \}$ such that
$\sum_i e_i \phi (f_i r) = r$ for each $r \in R$.
Now if $\phi(Rx) = 0$, then $x = \sum_i e_i \phi(f_i x) = 0$.
\end{proof}

We say that a Frobenius algebra
$A'$ is \textit{index one}
if  there is a Frobenius homomorphism $\phi: A \to K$  with
dual bases $\{ e_i \}$, $\{ f_i \}$ such that
$\sum_i e_i f_i = 1$. Then $\sum_i e_i \o f_i$ is a separability idempotent
and $A'$ is $K$-separable.
If $\phi$ is a trace, $A'$ is moreover strongly separable in Kanzaki's
sense \cite{KS}; e.g., finite separable field extensions and indeed
 separable commutative algebras are strongly separable, hence
index one Frobenius.  However, the $2 \times 2$ matrix
example in characteristic two in \cite[Footnote 3]{KN}, also given below, is
an index one Frobenius
algebra which is not strongly separable.  In fact, separable algebras
over fieldsare all index one Frobenius as we see next.

\begin{pro}
A $K$-separable algebra $R$ is index one Frobenius.
\end{pro}
\begin{proof}
More generally, a Frobenius extension $M' | N'$
is said to be of \textit{index one} if there is Frobenius homomorphism $E: M'
\to N'$ with Watatani index $[M': N']_E = \sum_i x_i y_i = 1_{M'}$,
where $\{ x_i \}$, $\{ y_i \}$ are arbitrary dual bases
for $E$.  This is quite easily checked to be a transitive notion
(cf.\ \cite[Prop. 2.1]{KN}).  It is also easy to check that, first,
 $K \times \cdots \times K$
is an index one Frobenius $K$-algebra via $$(\lambda_1,\ldots, \lambda_n)
\mapsto \sum_i \lambda_i,$$ and, second, that the tensor algebra
$A' \o_K B'$ is index one Frobeniusif $A'$ and $B'$ are so.

We recall the general fact of Hirata-Sugano that a Frobenius extension
$M' | N'$
with data $E,x_i,y_i$ as above is separable iff there is $d \in C_{M'}(N')$
suchthat $\sum_i x_idy_i=$ $1$ (cf.\ \cite{K}).
Next we consider a separable division algebra $D'$ (i.e., its center
$Z(D')$ is a  separable
field extension of $K$).  Certainly, $D'$ is a Frobenius algebra with
Frobeniushomorphism $\phi: D' \to K$ and dual bases $e_i, f_i$.
Then there is invertible
$d \in D'$ such that $\sum_i e_i df_i = 1$ where $e_i, df_i$ are dual bases of
index one forthe Frobenius homomorphism $\phi d^{-1}$.

Since $R$ is semisimple, we have the main Wedderburn theorem stating that
$$R \cong M_{n_1}(D_1) \times \cdots \times M_{n_t}(D_t),$$
where each $D_i$ is a separable division algebra.  Since $D_i$ is an index
one Frobeniusalgebra and $M_{n_i}(D_i) \cong M_{n_i}(K) \o D_i$, it remains
to show thateach full matrix algebra $M_n(K)$ is index one Frobenius, a
proof in three parts.

If the characteristic of $K$ is zero, or prime $p$ such that $p$ does
not divide$n$, then $M_n(K)$ is Kanzaki strongly separable \cite{KS} and
therefore index oneFrobenius.

If $\mathrm{char}\, K \geq 3$ and divides $n$, we modify the usual Frobenius
coordinates
$\phi(X) = \sum_i X_{ii}$ for $X \in M_n(K)$ with dual bases given by matrix
units
$e_{ij}, e_{ji}$ as follows.  Let $D$ be a diagonal matrix with first diagonal
entry $d_1 := 2$, the rest $d_i := 1$, then $\det D \neq 0$, and
we may consider the Frobenius
homomorphism $\phi D^{-1}$ with index
$$
\sum_{i,j} e_{ij}D e_{ji} = \sum_{i,j} d_j e_{ii} = (d_1 + \cdots + d_n) 1
                          = 1.
$$

Finally, if $\mathrm{char}\, K = 2$ and $n = 2^q m$ where $\gcd(m,2) = 1$ and
$q \geq 1$,
we note that $M_n(K) \cong M_2(K) \o M_{2^{q-1}m}(K)$.  The proof then
proceeds by induction
on $q$ if we show $M_2(K)$ to be index one Frobenius.  But this is the
content of
\cite[Footnote 3]{KN} where it is noted that
$\phi(X) = X_{11} + X_{12} + X_{21}$has dual bases given by
$$
e_{11} \otimes e_{21} + e_{12} \otimes e_{11} + e_{12} \otimes e_{21} + e_{22}
\otimes e_{12}+ e_{22} \otimes e_{22} + e_{21} \otimes e_{22},
$$
clearly of index one.
\end{proof}

For the rest of the section, $\phi, e_i, f_i$ will denote
  index one Frobenius coordinates for a separable algebra $R$.
Next we provide a converse to \cite{EN} and
a proof for a left-handed version of \cite[1.6]{Sz},
which states that a bialgebroid
over an index one Frobenius algebra is a weak bialgebra.

\begin{pro}\label{pro: WHA}
Suppose $R$ is a separable algebra and $(A,R,s,t,\Delta_A, \eps_A)$
is a bialgebroid in the category of $K$-algebras.
Then there is a weak bialgebra structure $(A, \Delta, \eps)$ defined
below.
\end{pro}
\begin{proof}
We denote $\Delta_A(a) = a\1 \o_R\, a\2$ as usual.    Then
we define the weak coproduct $\Delta: A \to A \o_K A$
 and weak counit $\eps: A \to K$ by
\begin{equation}
\Delta(a) := \sum_i t(e_i) a\1 \o s(f_i) a\2,
\end{equation}
\begin{equation}
\eps := \phi \eps_A.
\end{equation}

We first check that $(A, \Delta, \eps)$ is a coalgebra by applying
Eqs.~\ref{eq:del s} and~\ref{eq:del t}.
$$
(\id \o \cop)\cop(a) = \sum_{i,j} t(e_i) a\1 \o t(e_j)s(f_i)a\2 \o s(f_j)a\3
= (\cop \o \id)\cop(a),
$$
since $t(r)s(r') = s(r')t(r)$ for $r,r' \in R$. Next,
\begin{eqnarray*}
( \eps \o \id)\cop(a) &=& \sum_i \phi(\eps_A(t(e_i)a\1))s(f_i)a\2 = \\
\sum_i \phi(\eps_A(a\1)e_i)s(f_i)a\2 &=& s(\eps_A(a\1))a\2 = a,
\end{eqnarray*}
and similarly
$$ (\id \o \eps) \cop(a) = \sum_i t(e_i)a\1 \phi(f_i \eps_A(a\2)) = a.
$$

Next we use property (iii) for left
bialgebroids and $\sum_j e_j f_j = 1$ to show that $\cop$
is multiplicative: for $a,b \in A$,
\begin{eqnarray*}
 \cop(a) \cop(b) &=&
\sum_{i,j} t(e_i) a\1 t(e_j) b\1 \o s(f_i) a\2 s(f_j) b\2 \\
&= &  \sum_{i,j} t(e_i) a\1 b\1 \o s(f_i) a\2 s(e_j f_j) b\2 = \cop(ab).
\end{eqnarray*}

Note that $\cop(1) = \sum_i t(e_i) \o s(f_i)$. Then:
$$ (\cop(1) \o \id)(\id \o \cop(1)) = \sum_i t(e_i) \o s(f_i) t(e_j) \o s(f_j)
= (\id \o \cop(1))(\cop(1) \o \id) = \cop^2(1) $$
since $t(r)s(r') = s(r')t(r)$.

Next we establish the
weak multiplicativity of the counit with the help of property (vii):
\begin{eqnarray*}
\eps(ab\1)\eps(b\2 c) & = & \sum_i \phi(\eps_A(a t(e_i) b\1))
\phi(\eps_A(s(f_i)b\2c)) =\\
\phi( \eps_A(a t(\eps_A(b\2 c))b\1)) & = & \phi(\eps_A(a
t(\eps_A(b\1)\eps_A(b\2c)))) = \\
\phi(\eps_A(abc)) & = & \eps(abc)
\end{eqnarray*}
\begin{eqnarray*}
\eps(ab\2)\eps(b\1 c) & = & \sum_i \phi(\eps_A(a s(f_i) b\2))
\phi(\eps_A(t(e_i)b\1c)) =\\
\phi( \eps_A(a s(\eps_A(b\1 c))b\2)) & = &
\phi(\eps_A(a s(\eps_A(b\1c)\eps_A(b\2)))) = \\
\phi(\eps_A(a s(\eps_A(t(\eps_A(b\2))b\1 c)))) & = & \eps(abc)
\end{eqnarray*}
Thus $(A, \cop, \eps)$ is a weak bialgebra (cf.\ \cite{BNS}).
\end{proof}

The corresponding formulas for weak coproduct and weak counit for
a \textit{right} bialgebroid $(B,R,s,t,\cop_B,\eps_B)$ are given (as
in \cite[1.6]{Sz}) by
\begin{equation}
\cop(b) = \sum_i b\1 s(e_i) \o b\2 t(f_i), \ \ \ \eps = \phi \circ \eps_B.
\end{equation}

Let $K$ again denote a field.
\begin{thm}
\label{thm: A is WHA}
If $M |N$ is a D2 Frobenius extension of $K$-algebras  with centralizer
$R$ a separable algebra, then $A$ and $B$ are weak Hopf
algebras dual to one another.
\end{thm}
\begin{proof}
Again let $(\phi, e_j, f_j)$ denote an index one Frobenius coordinate
system for $R$.  By the proposition, $A$ and $B$ are weak bialgebras over
$K$ with weakcoproducts and weak counits given by: ($a \in A$, $b \in B$)
\begin{eqnarray}
\cop(a) &=& \sum_{i,j} a(-b^1_i)b^2_i e_j \o \lambda(f_j) \beta_i \\
\eps(a) &=& \phi(a(1)) \\
\cop(b)  & = & \sum_{i,j} (b^1_i \o_N b^2_i e_j) \o_K (f_j \beta_i(b^1) \o_N
b^2) \\
\eps(b) & = & \phi(b^1 b^2)
\end{eqnarray}

As we have seen in Lemma~\ref{lem: R-dual = K-dual} and Section~3,
there is a nondegenerate pairing obtained from composing
$B \cong \Hom_R(A,R) \cong \Hom_K(A,K)$ given by
(one of two choices):
\begin{equation}
\bra b,a \ket := \phi ( b^1 a(b^2)).
\end{equation}

We check that the weak bialgebra structures on $A$ and $B$
 are dual to one another with respect to $\bra , \ket$: ($b,c \in B$)

\begin{eqnarray*}
\bra b, a\1 \ket \bra c, a\2 \ket & = & \sum_{i,j}
\phi(b^1 a(b^2 b^1_j)b^2_j e_i) \phi(c^1 f_i \beta_j(c^2)) \\
& = & \phi(c^1 b^1 a(b^2 b^1_j)b^2_j \beta_j(c^2)) \\
& = & \phi(c^1 b^1 a(b^2 c^2)) = \bra bc, a \ket
\end{eqnarray*}
since $\sum_i \phi(re_i)f_i = r$ for
$r \in R$ and $bc = c^1 b^1 \o b^2 c^2$.  Also, ($a,\alpha \in A$)
\begin{eqnarray*}
\bra b\1, a \ket \bra b\2, \alpha \ket & = &
\phi(b^1_i a(b^2_i e_j)) \phi(f_j \beta_i(b^1) \alpha(b^2)) \\
& = & \phi(b^1_i a(b^2_i \beta_i(b^1) \alpha(b^2))) \\
& = & \phi(b^1 a \alpha(b^2)) = \bra b, a \alpha \ket
\end{eqnarray*}
Finally,  $$\bra 1_B, a \ket = \phi(a(1)) = \eps(a)$$
and $$\bra b, 1_A \ket = \phi(b^1 b^2) = \eps(b). $$
Hence, $A$ and $B$ are dual weak bialgebras.

We note that $A$ and $B$ are finite dimensional, for $R$ is
finite dimensional by assumption and $M| N$, therefore $M_1 | M$,
are finitely generated extensions.  But $M_1 | M$ has
dual bases identical (up to the isomorphism $\psi_A$)
with dual bases for the extension $A | R$ (cf.\ Section~6),
whence $A$ has finite $K$-dimension.
It then suffices to show that $A$ has an antipode.  For this we
will make use of the $\Leftarrow$ part of
 \cite[Theorem 4.1]{V}:
\begin{quote}
A finite dimensional weak bialgebra $A$ is a weak Hopf algebra iff
there is a nondegenerate left integral in $A$.
\end{quote}

We compute the projection $\Pi^L : A \to A^L$ onto the left (or target)
subalgebra of $A$:
\begin{eqnarray}
\Pi^L(a) = \eps(1\1 a) 1\2 &=& \sum_i \eps(\rho(e_i)a) \lambda(f_i)
\nonumber \\
& = & \phi ( a(1) e_i) \lambda(f_i) = \lambda(a(1) ) \label{eq: left map}
\end{eqnarray}

Recall that an element $\ell \in A$ is \textit{left integral} if $a \ell =
\Pi^L(a) \ell$
for every $a \in A$.  The Frobenius homomorphism $E \in \Hom_{N-N}(M,N)
\to \End_{N-N}(M)$ is such an element: ($m \in M$)
\begin{equation}\label{eq: left integral}
 aE(m) = a(1) E(m) \ \ \Rightarrow \ \ aE = \lambda(a(1))E = \Pi^L(a)E
\end{equation}
for each $a \in A$.

Recall that a left integral $\ell$ in $A$
 is \textit{nondegenerate} if the
two maps of $B \to A$ given by
\begin{eqnarray*}
b & \mapsto & \ell \leftharpoonup b = \bra b, \ell\1 \ket \ell\2 \\
b & \mapsto & b \rightharpoonup \ell = \ell\1 \bra b, \ell\2 \ket
\end{eqnarray*}
are isomorphisms.  However, by Lemma~\ref{lem: left iff right nondegen}
and finite dimensionality,
it will suffice to check that the simplest of these two maps is surjective.
We compute: ($b \in B$)
$$
E \leftharpoonup b = \sum_{i,j} \phi(b^1 E(b^2 b^1_i)b^2_ie_j)f_j \beta_i
= b^1 E(b^2 -) = E \ract b .$$

Given $\alpha \in A$, choose $b := \sum_j \alpha(x_j) \o y_j \in B$
and note that $E \ract b = \alpha$.  Thus, $E$ is nondegenerate left
integral in $A$, and $A$ is a weak Hopf algebra.  By duality for
weak Hopf algebras, $B$ is isomorphic to the dual weak Hopf algebra of $A$.
\end{proof}

Under the hypothesis that
$N \into M$ is a Frobenius D2 extension with separable centalizer $R$
such that $M_N$ is balanced,
we have the following equivalence of separability for quantum algebra
and proper algebra extension.
\begin{cor}
$M | N$ is split (resp., separable) if and only if $A$ (resp., $B$) is
$K$-separable.
\end{cor}
\begin{proof}
Since $M | N$ is Frobenius, we have noted in Section~8 that $M | N$ is split
iff there is $d \in C_M(N)$ such that $E(d) = 1$.  But the element
of $A$ given by $Ed: x \mapsto
E(dx)$ ($x \in M$) is a left integral via a simple calculation
as in Eq.\ (\ref{eq: left integral}).  And by Eq.\ (\ref{eq: left map})
$Ed$ is a normalized left integral, whose existence is
 by \cite[Theorem 3.13]{BNS} equivalent to
$A$ being $K$-separable.

Conversely, if $\ell \in A$ is a normalized left integral, the
mapping $m \mapsto \ell \lact m$ is easily seen to induce an
$N$-bimodule splitting
map for $N \into M$ since $N = M^A$. Whence $M | N$ is a split extension.

If $M | N$ is separable Frobenius (though not necessarily proper), then
there is $d \in C_M(N)$ such that $\sum_j x_j d y_j = 1$.  We
compute that $\sum_j x_j \o dy_j$ is a normalized right integral in $B$
by noting that
$$
\Pi^R(b) = 1\1 \eps(b 1\2) = \sum_{i,j} b^1_i \o b^2_i e_j \phi(f_j
\beta_i(1)b^1 b^2) = 1_B \eps_R(b),
$$
and a computation as in Eq.\ (\ref{eq: right integral}).

Conversely, given a normalized right integral $\beta \in B$,
we induce a splitting map $\cdot \ract \beta$ for
$\rho: M \to M_1^{\rm op}$.  Then $M_1 | M$ is split.  Since $C_M(N)$
is anti-isomorphic with $C_{M_1}(M)$ via $$d \mapsto
\sum_i x_i de_1 y_i,$$
there is $d' := \sum_i x_i de_1 y_i \in C_{M_1}(M)$ such that
$E_M(d') = \sum_i x_i d y_i = 1$ where $d \in C_M(N)$,
i.e., $M | N$ is separable.
\end{proof}

\end{document}